\documentclass[10pt]{article}
\usepackage{amsmath}
\usepackage{latexsym}
\usepackage{amssymb}
\usepackage{amsfonts}
\usepackage{amsthm}
\usepackage{tikz}
\newcommand*\circled[1]{\tikz[baseline=(char.base)]{\node[shape=circle,draw,inner sep=2pt](char){#1};}}
\title{RELATIVE  PARAMETRIZATION  OF  \\  LINEAR  MULTIDIMENSIONAL  SYSTEMS}
\author{J.F. Pommaret \\ CERMICS, Ecole Nationale des Ponts et
  Chauss\'ees,\\6/8 Av. Blaise Pascal, 77455 Marne-la-Vall\'ee Cedex 02,
  France \\
E-mail: jean-francois.pommaret@wanadoo.fr, pommaret@cermics.enpc.fr \\
URL: http://cermics.enpc.fr/${\sim}$pommaret/home.html }
\date{  }
\begin{document}
\maketitle

\noindent
{\bf  ABSTRACT} :  \\

\noindent
\hspace*{4mm} In the last chapter of his book "{\it The Algebraic Theory of Modular Systems} " published in 1916, F. S. Macaulay developped specific techniques for dealing with "{\it unmixed polynomial ideals} " by introducing what he called "{\it inverse systems} ". The purpose of this paper is to extend such a point of view to differential modules defined by linear multidimensional systems, that is by linear systems of ordinary differential (OD) or partial differential (PD) equations of any order, with any number of independent variables, any number of unknowns and even with variable coefficients.  \\
    The {\it first and main idea} is to replace unmixed polynomial ideals by {\it pure differential modules}.  \\
    The {\it second idea} is to notice that a module is $0$-pure if and only if it is {\it torsion-free} and thus if and only if it admits an " {\it absolute parametrization} " by means of arbitrary potential like functions, or, equivalently, if it can be embedded into a free module by means of an " {\it absolute localization} ".  \\
     The {\it third idea} is to refer to a difficult theorem of algebraic analysis saying that an $r$-pure module can be embedded into a module of projective dimension $r$, that is a module admitting a projective resolution with exactly $r$ operators.  \\
     The {\it fourth and final idea} is to establish a link between the use of extension modules for such a purpose and specific formal properties of the underlying multidimensional system through the use of involution and a "{\it relative localization} " leading to a " {\it relative parametrization} ".  \\
     The paper is written in a rather effective self-contained way and we provide many explicit examples that should become test examples for a future use of computer algebra.  \\
     
{\bf  KEY WORDS} :

\noindent
Unmixed ideal, algebraic analysis, homological algebra, extension module, projective dimension, torsion-free module, pure module, characteristic variety, formal integrability, involution, Spencer operator, inverse system.  \\

\noindent
{\bf 1)  INTRODUCTION} :   \\

Let $D=K[d_1,...,d_n]=K[d]$ be the ring of differential operators with coefficients in a differential field $K$ with $n$ commuting derivations ${\partial}_1,...,{\partial}_n$ and commutation relations $d_ia=ad_i+{\partial}_ia,\forall a\in K$. If $y^1,...,y^m$ are $m$ differential indeterminates, we may identify $Dy^1+...+Dy^m=Dy$ with $D^m$ and consider the finitely presented left differential module $M$ with presentation $D^p\rightarrow D^m\rightarrow M \rightarrow 0$ determined by a given linear multidimensional system with $n$ independent variables, $m$ unknowns and $p$ equations. Applying the functor $hom_D(\bullet,D)$, we get the exact sequence $0\rightarrow hom_D(M,D)\rightarrow D^m\rightarrow D^p \longrightarrow N \longrightarrow 0$ of {\it right differential modules} that can be transformed by a side-changing functor to an exact sequence of finitely generated {\it left differential modules}. This new presentation corresponds to the {\it formal adjoint} $ad({\cal{D}})$ of the linear differential operator $\cal{D}$ determined by the initial presentation but now with $p$ unknowns and $m$ equations, obtaining therefore a new finitely generated {\it left differential module} $N$ and we may consider $hom_D(M,D)$ as the {\it module of equations} of the {\it compatibility conditions} (CC) of $ad({\cal{D}})$, a result not evident at all (Compare to [24]). Using now a maximum free submodule $0 \longrightarrow D^l \longrightarrow hom_D(M,D)$ and repeating this standard procedure while using the well known fact that $ad(ad({\cal{D}}))={\cal{D}}$, we obtain therefore an embedding $0\rightarrow hom_D(hom_D(M,D),D)\rightarrow D^l$ of left differential modules for a certain integer $1\leq l<m$ because $K$ is a field and thus $D$ is a noetherian bimodule over itself, a result leading to $l=rk_D(hom_D(M,D))=rk_D(M)< m$ as in ([19], p 178,201)(See section 3 for the definition of the {\it differential rank}$rk_D$). Now, the kernel of the map $\epsilon:M\rightarrow hom_D(hom_D(M,D),D):m\rightarrow \epsilon(m)(f)= f(m),\forall f\in hom_D(M,D)$ is the torsion submodule $t(M)\subseteq M$ and $\epsilon$ is injective if and only if $M$ is torsion-free, that is $t(M)=0$. In that case, we obtain by composition an embedding $0\rightarrow M \rightarrow D^l$ of $M$ into a free module that can also be obtained by localization if we introduce the ring of fractions $S^{-1}D=DS^{-1}$ when $S=D-\{0\}$. This result is quite important for applications as it provides a (minimum) parametrization of the linear differential operator $\cal{D}$ and amounts to the controllability of a classical control system when $n=1$ ([24], p 258). This parametrization will be called an "{\it absolute parametrization} " as it only involves arbitrary "{\it potential-like} " functions (See [1], [9], [18], [19], [20], [24], [25] and [32] for more details and examples, in particular that of Einstein equations).   \\

The purpose of this paper is to extend suh a result to a much more general situation, that is {\it when} $M$ {\it is not torsion-free}, by using unexpected results first found by F.S. Macaulay in 1916 through his study of "{\it inverse systems} " for "{\it unmixed polynomial ideals} ".  \\

For this we define the {\it purity filtration} :\\
    \[     0=t_n(M) \subseteq t_{n-1}(M) \subseteq  ... \subseteq t_1(M) \subseteq t_0(M)=t(M)\subseteq M  \]
by introducing $t_r(M)=\{m\in M\mid cd(Dm)>r\}$ where the codimension of $Dm$ is $n$ minus the dimension of the characteristic variety determined by $m$ in the corresponding system for one unknown. The module $M$ is said to be $r$-{\it pure} if $t_r(M)=0, t_{r-1}(M)=M$ or, equivalently, if $cd(M)=cd(N)=r, \forall N\subset M$ and a torsion-free module is a 0-pure module. Moreover, when $K=k=cst(K)$ is a field of constants and $m=1$, a pure module is {\it unmixed} in the sense of Macaulay, that is defined by an ideal having an equidimensional primary decomposition. \\

\noindent
{\bf Example 1.1} :  As an elementary example with $K=k=\mathbb{Q}, m=1,n=2, p=2$, the differential module defined by $d_{22}y=0,d_{12}y=0$ is not pure because $z'=d_2y$ satisfies $d_2z'=0,d_1z'=0$ while $z"=d_1y$ only satisfies $d_2z"=0$ and $(({\chi}_2)^2,{\chi}_1{\chi}_2)=({\chi}_1)\cap ({\chi}_1,{\chi}_2)^2$. We obtain therefore the purity filtration $ 0 = t_2(M) \subset  t_1(M) \subset t_0(M)=t(M)=M $ with strict inclusions as $ 0\neq z' \in t_1(M) $ while $z" \in t_0(M)$ but $z" \notin t_1(M) $.   \\

From the few (difficult) references ([1],[9],[15],[18],[27]) dealing with extension modules $ext^r(M)=ext^r_D(M,D)$ and purity in the framework of algebraic analysis, it is known that $M$ is $r$-pure if and only if there is an embedding $0 \rightarrow M \rightarrow ext^r_D(ext^r_D(M,D),D)$. Indeed, the case $r=0$ 
is exactly the one already considered because $ext^0_D(M,D)=hom_D(M,D)$ and the ker/coker exact sequence:\\
\[   0 \longrightarrow ext^1(N) \longrightarrow M \longrightarrow ext^0(ext^0(M)) \longrightarrow ext^2(N) \longrightarrow 0  \]
allows to test the torsion-free property of $M$ in actual practice by using the double-duality formula $t(M)=ext^1(N)$ as in ([19]). Also, when $r\geq 1$, a similar construction that we shall recall and illustrate in section 4 provides a finitely generated module $L$ with {\it projective dimension} $pd_D(L)=r$, that is a minimum resolution of $L$ with only $r$ operators, and an embedding $0 \rightarrow M \rightarrow L $ that allows to exhibit a {\it relative parametrization} of $\cal{D}$ because now the parametrizing potential-like functions are no longer arbitrary but must only depend on arbitrary functions of $n-r$ variables. \\

\noindent
{\bf Example 1.2} : With $K=k=\mathbb{Q},m=2, n=3, r=1$, the differential module $M$ defined by the involutive system ${\Phi}^1\equiv d_3y^1=0, {\Phi}^2\equiv d_3y^2=0, {\Phi}^3\equiv d_2y^1-d_1y^2=0$ is $1$-pure and admits the resolution $ 0 \longrightarrow D \longrightarrow D^3 \longrightarrow D^2 \longrightarrow M \longrightarrow 0 $. The differential module $L$ defined by the system $d_3z=0$ is also $1$-pure and admits the resolution $ 0 \longrightarrow D \longrightarrow D \longrightarrow L \longrightarrow 0 $. We finally obtain the relative parametrization $y^1=d_1z, y^2=d_2z$ providing the strict inclusion $M\subset L$.\\

In a simple way, this result can be considered as a measure of how far a module is from being projective, recalling that a module $P$ is 
{\it projective} if there exists another (projective) module $Q$ and a free module $F$ such that $P\oplus Q\simeq F$.\\

We adapt the "{\it relative localization}" technique used by Macaulay and combine it with the "{\it involution}" technique used in the formal theory of systems of partial differential equations in order to obtain an explicit procedure for determining $L$ when $M$ is given. Many examples will illustrate these new methods that avoid the previous abstract arguments based on "{\it double duality}". In particular, original non-commutative examples will also be presented. However, we point out the fact that the latter method can be adapted without any change to the case of systems with variables coefficients as it only depends on the use of adjoint operators but the following example will explain by itself the type of difficulty involved.  \\

\noindent
{\bf Example 1.3} :  Starting now with $K=\mathbb{Q} (x^1,x^2), m=2, n=3, r=1$, the new differential module $M$ defined by  $d_3y^1=0, d_3y^2=0, d_2y^1-d_1y^2+x^2y^2=0$ is also $1$-pure and the differential module $L$ is again defined by $d_3z=0$ as in the previous example. However we obtain the 
{\it totally different} relative parametrization $y^1=d_{12}z-x^2d_2z+z, y^2=d_{22}z$ providing the strict inclusion $M\subset L$. More generally, we may consider a constant parameter $a\in k=\mathbb{Q}$ and consider the new system $d_3y^1=0, d_3y^2=0, d_2y^1-d_1y^2+ax^2y^2=0$ depending on $a$. For $a=0$ we find back the case of the previous example and we let the reader wonder why the situation only changes when $a\neq 0$. \\

The content of the paper is just following the introduction.  \\

In section 2 we recall the definitions and results from the formal theory of systems of OD/PD equations that will be crucially used in the sequel. We pay a particular emphasize to the definition of involution and the way to introduce the Spencer operator in this framework. We also study the possibility and 
difficulty to use computer algebra in this framework.\\

In section 3 we recall the basic tools needed from module theory and homological algebra in a way adapted to our purpose, in particular the definition of the extension modules, and provide a few of their properties which, though well known by specialists of algebraic analysis, cannot be found easily in the literature. Meanwhile, we provide a few links with the preceding section which are not so well known. Many explicit examples will illustrate the main concepts in the commutative (constant field $k$) and the non-commutative (differential field $K$) framework.  \\

In section 4 we shall recall the proof of the theorem already quoted showing how to embed an $r$-pure module $M$ into another module $L$ with projrective dimension equal to $r$. We shall provide for the first time explicit computations of this result in order to point out the difficulty encountered in such a procedure as a motivation for avoiding it.   \\

In section 5 we extend the work of Macaulay, showing why {\it only pure modules can fit with relative localization} in a coherent way with what happens for torsion-free modules. Meanwhile, we shall extend for the first time this work to the non-commutative framework, showing in particular that the operator introduced by Macaulay ([11], \S 60) for studying inverse systems is nothing else than the Spencer operator. Many explicit examples, including highly non-trivial ones provided by Macaulay himself, will be fully treated in such a way that any engineer, even with a poor knowledge of homological algebra, will nevertheless become intuitively able to understand and apply these new techniques without reading the previous sections, just comparing to the way the same examples have been treated in section 4 by means of another approach.  \\
  
\noindent
{\bf 2)  TOOLS FROM  SYSTEM THEORY} :  \\

If $X$ is a manifold of dimension $n$ with local coordinates $(x)=(x^1, ... x^n)$, we denote as usual by $T=T(X)$ the {\it tangent bundle} of $X$, by $T^*=T^*(X)$ the {\it cotangent bundle}, by ${\wedge}^rT^*$ the {\it bundle of r-forms} and by $S_qT^*$ the {\it bundle of q-symmetric tensors}. More generally, let $E$ be a {\it vector bundle} over $X$, that is (roughly) a manifold with local coordinates $(x^i,y^k)$ for $i=1,...,n$ and $k=1,...,m$ simply denoted by $(x,y)$, {\it projection} $\pi:E\rightarrow X:(x,y)\rightarrow (x)$ and changes of local coordinates $\bar{x}=\varphi(x), \bar{y}=A(x)y$. If $E$ and $F$ are two vector 
bundles  over $X$ with respective local coordinates $(x,y)$ and $(x,z)$, we denote by $E{\times}_X F$ the {\it fibered product} of $E$ and $F$ over $X$ as the new vector bundle over $X$ with local coordinates $(x,y,z)$. We denote by $f:X\rightarrow E: (x)\rightarrow (x,y=f(x))$ a global {\it section} of $E$, that is a map such that $\pi\circ f=id_X$ but local sections over an open set $U\subset X$ may also be considered when needed. Under a change of coordinates, a section transforms like $\bar{f}(\varphi(x))=A(x)f(x)$ and the derivatives transform like:\\
\[   \frac{\partial{\bar{f}}^l}{\partial{\bar{x}}^r}(\varphi(x)){\partial}_i{\varphi}^r(x)=({\partial}_iA^l_k(x))f^k(x)+A^l_k(x){\partial}_if^k(x)  \]
We may introduce new coordinates $(x^i,y^k,y^k_i)$ transforming like:\\
\[ {\bar{y}}^l_r{\partial}_i{\varphi}^r(x)=({\partial}_iA^l_k(x))y^k+A^l_k(x)y^k_i  \]
We shall denote by $J_q(E)$ the {\it q-jet bundle} of $E$ with local coordinates $(x^i, y^k, y^k_i, y^k_{ij},...)=(x,y_q)$ called {\it jet coordinates} and sections $f_q:(x)\rightarrow (x,f^k(x), f^k_i(x), f^k_{ij}(x), ...)=(x,f_q(x))$ transforming like the sections $j_q(f):(x) \rightarrow (x,f^k(x), {\partial}_if^k(x), {\partial}_{ij}f^k(x), ...)=(x,j_q(f)(x))$ where both $f_q$ and $j_q(f)$ are over the section $f$ of $E$. Of course $J_q(E)$ is a vector bundle over $X$ with projection ${\pi}_q$ while $J_{q+r}(E)$ is a vector bundle over $J_q(E)$ with projection ${\pi}^{q+r}_q, \forall r\geq 0$.\\

\noindent
{\bf DEFINITION 2.1}: A {\it linear system} of order $q$ on $E$ is a vector sub-bundle $R_q\subset J_q(E)$ and a {\it solution} of $R_q$ is a section $f$ of $E$ such that $j_q(f)$ is a section of $R_q$.\\

Let $\mu=({\mu}_1,...,{\mu}_n)$ be a multi-index with {\it length} ${\mid}\mu{\mid}={\mu}_1+...+{\mu}_n$, {\it class} $i$ if ${\mu}_1=...={\mu}_{i-1}=0,{\mu}_i\neq 0$ and $\mu +1_i=({\mu}_1,...,{\mu}_{i-1},{\mu}_i +1, {\mu}_{i+1},...,{\mu}_n)$. We set $y_q=\{y^k_{\mu}{\mid} 1\leq k\leq m, 0\leq {\mid}\mu{\mid}\leq q\}$ with $y^k_{\mu}=y^k$ when ${\mid}\mu{\mid}=0$. If $E$ is a vector bundle over $X$ with local coordinates $(x^i,y^k)$ for $i=1,...,n$ and $k=1,...,m$, we denote by $J_q(E)$ the $q$-{\it jet bundle} of $E$ with local coordinates simply denoted by $(x,y_q)$ and {\it sections} $f_q:(x)\rightarrow (x,f^k(x), f^k_i(x), f^k_{ij}(x),...)$ transforming like the section $j_q(f):(x)\rightarrow (x,f^k(x),{\partial}_if^k(x),{\partial}_{ij}f^k(x),...)$ when $f$ is an arbitrary section of $E$. Then both $f_q\in J_q(E)$ and $j_q(f)\in J_q(E)$ are over $f\in E$ and the {\it Spencer operator} just allows to distinguish them by introducing a kind of "{\it difference}" through the operator $D:J_{q+1}(E)\rightarrow T^*\otimes J_q(E): f_{q+1}\rightarrow j_1(f_q)-f_{q+1}$ with local components $({\partial}_if^k(x)-f^k_i(x), {\partial}_if^k_j(x)-f^k_{ij}(x),...) $ and more generally $(Df_{q+1})^k_{\mu,i}(x)={\partial}_if^k_{\mu}(x)-f^k_{\mu+1_i}(x)$. In a symbolic way, {\it when changes of coordinates are not involved}, it is sometimes useful to write down the components of $D$ in the form $d_i={\partial}_i-{\delta}_i$ and the restriction of $D$ to the kernel $S_{q+1}T^*\otimes E$ of the canonical projection ${\pi}^{q+1}_q:J_{q+1}(E)\rightarrow J_q(E)$ is {\it minus} the {\it Spencer map} $\delta=dx^i\wedge {\delta}_i:S_{q+1}T^*\otimes E\rightarrow T^*\otimes S_qT^*\otimes E$. The kernel of $D$ is made by sections such that $f_{q+1}=j_1(f_q)=j_2(f_{q-1})=...=j_{q+1}(f)$. Finally, if $R_q\subset J_q(E)$ is a {\it system} of order $q$ on $E$ locally defined by linear equations ${\Phi}^{\tau}(x,y_q)\equiv a^{\tau\mu}_k(x)y^k_{\mu}=0$ and local coordinates $(x,z)$ for the parametric jets up to order $q$, the $r$-{\it prolongation} $R_{q+r}={\rho}_r(R_q)=J_r(R_q)\cap J_{q+r}(E)\subset J_r(J_q(E))$ is locally defined when $r=1$ by the linear equations ${\Phi}^{\tau}(x,y_q)=0, d_i{\Phi}^{\tau}(x,y_{q+1})\equiv a^{\tau\mu}_k(x)y^k_{\mu+1_i}+{\partial}_ia^{\tau\mu}_k(x)y^k_{\mu}=0$ and has {\it symbol} $g_{q+r}=R_{q+r}\cap S_{q+r}T^*\otimes E\subset J_{q+r}(E)$ if one looks at the {\it top order terms}. If $f_{q+1}\in R_{q+1}$ is over $f_q\in R_q$, differentiating the identity $a^{\tau\mu}_k(x)f^k_{\mu}(x)\equiv 0$ with respect to $x^i$ and substracting the identity $a^{\tau\mu}_k(x)f^k_{\mu+1_i}(x)+{\partial}_ia^{\tau\mu}_k(x)f^k_{\mu}(x)\equiv 0$, we obtain the identity $a^{\tau\mu}_k(x)({\partial}_if^k_{\mu}(x)-f^k_{\mu+1_i}(x))\equiv 0$ and thus the restriction $D:R_{q+1}\rightarrow T^*\otimes R_q$ ([17],[18],[30]). \\
    
\noindent
{\bf DEFINITION 2.2}: $R_q$ is said to be {\it formally integrable} when the restriction ${\pi}^{q+r+1}_{q+r}:R_{q+r+1}\rightarrow R_{q+r} $ is an epimorphism $\forall r\geq 0$ or, equivalently, when all the equations of order $q+r$ are obtained by $r$ prolongations only $\forall r\geq 0$. In that case, $R_{q+1}\subset J_1(R_q)$ is a canonical equivalent formally integrable first order system on $R_q$ {\it with no zero order equations}, called the {\it Spencer form}.\\

Finding an intrinsic test has been achieved by D.C. Spencer in 1970 ([30]) along coordinate dependent lines sketched by M. Janet in 1920 ([7]) and W. Gr\"{o}bner in 1940 ([4],[6]). The key ingredient, missing in the old approach, is provided by the following definition.\\

Let $T^*$ be the cotangent vector bundle of 1-forms on $X$ and ${\wedge}^sT^*$ be the vector bundle of s-forms on $X$ with usual bases $\{dx^I=dx^{i_1}\wedge ... \wedge dx^{i_s}\}$ where we have set $I=(i_1< ... <i_s)$. Moreover, introducing the {\it exterior derivative} $d: {\wedge}^sT^* \longrightarrow {\wedge}^{s+1}T^*:\omega={\omega}_I(x)dx^I \longrightarrow d\omega={\partial}_i{\omega}_I(x)dx^i\wedge dx^I$, we have $d^2=d\circ d=0$ and may introduce the {\it Poincar\'{e} sequence}:        \\

\noindent
\[   {\wedge}^0T^* \stackrel{d}{\longrightarrow} {\wedge}^1T^* \stackrel{d}{\longrightarrow} {\wedge}^2T^* \stackrel{d}{\longrightarrow} ... \stackrel{d}{\longrightarrow}{\wedge}^nT^* \longrightarrow 0  \]

\noindent
{\bf PROPOSITION 2.3}: There exists a map $\delta:{\wedge}^sT^*\otimes S_{q+1}T^*\otimes E\rightarrow {\wedge}^{s+1}T^*\otimes S_qT^*\otimes E$ which restricts to $\delta:{\wedge}^sT^*\otimes g_{q+1}\rightarrow {\wedge}^{s+1}T^*\otimes g_q$ and ${\delta}^2=\delta\circ\delta=0$.\\

{\it Proof}: Let us introduce the family of s-forms $\omega=\{ {\omega}^k_{\mu}=v^k_{\mu,i}dx^I\}$ and set $(\delta\omega)^k_{\mu}=dx^i\wedge{\omega}^k_{\mu+1_i}$. We obtain at once $({\delta}^2\omega)^k_{\mu}=dx^i\wedge dx^j\wedge{\omega}^k_{\mu+1_i+1_j}=0$. \\ 
\hspace*{12cm}  Q.E.D.    \\

The kernel of each $\delta$ in the first case is equal to the image of the preceding $\delta$ but this may no longer be true in the restricted case and we set:\\

\noindent
{\bf DEFINITION 2.4}: We denote by $H^s_{q+r}(g_q)$ the cohomology at ${\wedge}^sT^*\otimes g_{q+r}$ of the restricted $\delta$-sequence which only depends on $g_q$. The symbol $g_q$ is said to be s-{\it acyclic} if $H^1_{q+r}=...=H^s_{q+r}=0, \forall r\geq 0$, involutive if it is n-acyclic and {\it finite type} if $g_{q+r}=0$ becomes trivially involutive for r large enough.\\

\noindent
{\bf DEFINITION 2.5}: $R_q$ is said to be {\it involutive} when it is formally integrable and its symbol $g_q$ is involutive, that is to say all the sequences $... \stackrel{\delta}{\rightarrow} {\wedge}^sT^*\otimes g_{q+r}\stackrel{\delta}{\rightarrow}...$ are exact $\forall 0\leq s\leq n, \forall r\geq 0$. \\

Equivalently, the following procedure, {\it where one may have to change linearly the independent variables if necessary}, is the heart towards the next effective definition of involution. It is intrinsic even though it must be checked in a particular coordinate system called $\delta$-{\it regular} ([17],[18],[29]) and is particularly simple for first order systems without zero order equations.    \\

\noindent
$\bullet$ {\it Equations of class} $n$: Solve the maximum number ${\beta}^n_q$ of equations with respect to the jets of order $q$ and class $n$. Then call $(x^1,...,x^n)$ {\it multiplicative variables}.\\
\[  - - - - - - - - - - - - - - - -  \]
$\bullet$ {\it Equations of class} $i\geq 1$: Solve the maximum number ${\beta}^i_q$ of {\it remaining} equations with respect to the jets of order $q$ and class $i$. Then call $(x^1,...,x^i)$ {\it multiplicative variables} and $(x^{i+1},...,x^n)$ {\it non-multiplicative variables}.\\
\[ - - - - - - - - - - - - - - - - - \]
$\bullet$ {\it Remaining equations equations of order} $\leq q-1$: Call $(x^1,...,x^n)$ {\it non-multiplicative variables}.\\

\noindent
In actual practice, we shall use a {\it multiplicative board} where the multiplicative "variables" are represented by their index in upper left position while the non-multiplicative variables are represented by dots in lower right position.  \\

\noindent
{\bf DEFINITION 2.6}: A system of PD equations is said to be {\it involutive} if its first prolongation can be achieved by prolonging its equations only with respect to the corresponding multiplicative variables. In that case, we may introduce the {\it characters} ${\alpha}^i_q=m\frac{(q+n-i-1)!}{(q-1)!((n-i)!}-{\beta}^i_q$ for $i=1, ..., n$ and we have $dim(g_{q+1})={\alpha}^1_q+...+{\alpha}^n_q$. Moreover, one can exhibit the {\it Hilbert polynomial} $dim(R_{q+r})$ in $r$ with leading term $(\alpha/d!)r^d$ with $d\leq n$ when $\alpha$ is the smallest non-zero character in the case of an involutive symbol. Such a prolongation allows to compute {\it in a unique way} the principal ($pri$) jets from the parametric ($par$) other ones. This definition may also be applied to nonlinear systems as well.  \\

\noindent
{\bf REMARK 2.7}: For an involutive system with $\beta ={\beta}^n_q< m$, then $(y^{\beta +1},...,y^m)$ can be given arbitrarily and may constitute the {\it input} variables in control theory, though it is not necessary to make such a choice. {\it In this case}, the intrinsic number $\alpha={\alpha}^n_q=m-\beta> 0$ is called the $n$-{\it character} and is the system counterpart of the so-called "{\it differential transcendence degree}" in differential algebra. As we shall see in the next section, {\it the  smallest non-zero character and the number of zero characters are intrinsic numbers that cannot be known without bringing the system to involution} and we have 
${\alpha}^1_q\geq ... \geq {\alpha}^n_q\geq 0$. \\

\noindent
{\bf EXAMPLE 2.8}: ([11], \S 38, p 40 where one can find the first intuition of formal integrability) The primary ideal $\mathfrak{q}=(({\chi}_1)^2, {\chi}_1{\chi}_3-{\chi}_2)$ provides the system $y_{11}=0, y_{13}-y_2=0$ which is neither formally integrable nor involutive. Indeed, we get $d_3y_{11}-d_1(y_{13}-y_2)=y_{12}$ and $d_3y_{12}-d_2(y_{13}-y_2)=y_{22}$, that is to say {\it each first and second} prolongation does bring a new second order PD equation. Considering the new system $y_{22}=0, y_{12}=0, y_{13}-y_2=0, y_{11}=0$, the question is to decide whether this system is involutive or not. One could use Janet or Gr\"{o}bner algorithm but with no insight towards involution. In such a simple situation, as there is no PD equation of class $3$, two evident permutations of coordinates $(1,2,3)\rightarrow (3,2,1)$ or $(1,2,3) \rightarrow (2,3,1)$ both provide one equation of class $3$, $2$ equations of class $2$ and $1$ equation of clas $1$. It is then easy to check directly that the first permutation brings the involutive system $y_{33}=0, y_{23}=0, y_{22}=0, y_{13}-y_2=0$ that will be used in the sequel and we have ${\alpha}^3_2=0, {\alpha}^2_2=0, {\alpha}^1_2=2$.  \\

\noindent
{\bf EXAMPLE 2.9}: With $n=4, m=1, q=1, K=\mathbb{Q}(x^1,x^2,x^3,x^4)$, let us consider the system $R_1$:\\
\[    \hspace{3cm} \{  y_4-x^3y_2-y=0, \hspace{1cm} y_3-x^4y_1=0  \]
Again, the reader will check easily that the subsystem $R'_1\subset R_1$:\\
\[  \left\{  \begin{array}{rl}
u  \equiv &   y_4-x^3y_1-y=0  \\
v   \equiv  &  y_3-x^4y_1=0 \\
w  \equiv &    y_2-y_1=0 
\end{array}
\right.  \fbox{  $ \begin{array}{llll}
1 & 2 & 3 & 4 \\
1 & 2  & 3  & \bullet \\
1 & 2  & \bullet & \bullet 
\end{array}  $ }  \] 
namely the projection $R^{(1)}_1$ of $R_2$ to $R_1$, is formally integrable and even involutive with one equation of class 4, one equation of class 3 and one equation of class 2.\\

In the situation of the last remark, the following theorem will generalizing for PD control systems the well known first order Kalman form of OD control systems where the derivatives of the input do not appear ([27], VI,1.14, p 802). For this, {\it we just need to modify the Spencer form} and we provide the procedure that must be followed in the case of a first order involutive system with no zero order equation, for example an involutive Spencer form.\\

\noindent 
$\bullet$  Look at the equations of class $n$ solved with respect to $y^1_n,...,y^{\beta}_n$.\\
$\bullet$  Use integrations by part like:\\
\[ y^1_n-a(x)y^{\beta +1}_n=d_n(y^1-a(x)y^{\beta +1})+{\partial}_na(x)y^{\beta +1}={\bar{y}}^1_n+{\partial}_na(x)y^{\beta +1}  \]
$\bullet$  Modify $y^1,...,y^{\beta} $ to ${\bar{y}}^1,...,{\bar{y}}^{\beta}$ in order to "{\it absorb}" the various $y^{\beta +1}_n,...,y^m_n$ {\it only appearing in the equations of class} n.\\

We have the following unexpected result providing what we shall call {\it reduced Spencer form}:\\

\noindent
{\bf THEOREM 2.10}: The new equations of class n only contain $y^{\beta +1}_i,...,y^m_i$ with $0\leq i\leq n-1$ while the equations of class $1,...,n-1$ no longer contain $y^{\beta+1},...,y^m$ and their jets. Accordingly, as we shall see in the next section, any torsion element, if it exists, only depends on ${\bar{y}}^1,...,{\bar{y}}^{\beta}$.\\

{\it Proof}: The first assertion comes from the absorption procedure. Now, if $y^m$ or $y^m_i$ should appear in an equation of class $\leq n-1$, prolonging this equation with respect to the non-multiplicative variable $x^n$ should bring $y^m_n$ or $y^m_{in}$ and (here involution is essential) we should get a linear combination of equations of various classes prolonged with respect to $x^1,...,x^{n-1}$ {\it only}, but this is impossible and we get the desired reduced form.        \\
\hspace*{12cm} Q.E.D.  \\
When $R_q$ is involutive, the linear differential operator ${\cal{D}}:E\stackrel{j_q}{\rightarrow} J_q(E)\stackrel{\Phi}{\rightarrow} J_q(E)/R_q=F_0$ of order $q$ with space of solutions $\Theta\subset E$ is said to be {\it involutive} and one has the canonical {\it linear Janet sequence} ([17], p 144):\\
\[  0 \longrightarrow  \Theta \longrightarrow T \stackrel{\cal{D}}{\longrightarrow} F_0 \stackrel{{\cal{D}}_1}{\longrightarrow}F_1 \stackrel{{\cal{D}}_2}{\longrightarrow} ... \stackrel{{\cal{D}}_n}{\longrightarrow} F_n \longrightarrow 0   \]
where each other operator is first order involutive and generates the {\it compatibility conditions} (CC) of the preceding one. As the Janet sequence can be cut at any place, {\it the numbering of the Janet bundles has nothing to do with that of the Poincar\'{e} sequence} for the exterior derivative, contrary to what many physicists  believe. Moreover, the dimensions of the Janet bundles can be computed at once inductively from the board of multiplicative and non-multiplicative variables that can be exhibited for $\cal{D}$ by working out the board for ${\cal{D}}_1$ and so on. For this, the number of rows of this new board is the number of dots appearing in the initial board while the number $nb(i)$ of dots in the column $i$ just indicates the number of CC of class $i$ for $i=1, ... ,n$ with 
$nb(i) < nb(j), \forall i<j$. It follows that the successive first order operators ${\cal{D}}_1, ... , {\cal{D}}_n$ are {\it automatically} in reduced Spencer form. \\

\noindent
{\bf EXAMPLE 2.11}: Coming back to Example 2.9 and changing slightly our usual notations, we get for ${\cal{D}}_1$ the following first order involutive system 
of CC in reduced Spencer form:  \\

\[   \left\{  \begin{array}{ll}
{\phi}^3\equiv & d_4v-d_3u+x^4d_1u-x^3d_1v-v =0  \\
{\phi}^2\equiv & d_4w-d_2u+d_1u-x^3d_1w-w =0 \\
{\phi}^1 \equiv & d_3w-d_2v+d_1v-x^4d_1w =0
\end{array}
\right. \fbox{ $ \begin{array}{llll}
1 & 2 & 3 & 4  \\
1 & 2 & 3 &   4 \\
1 & 2 & 3 & \bullet
\end{array} $ } \]

\noindent
as $d_4u$ does not appear in ${\phi}^2$ and ${\phi}^3$ while $u$ does not appear in ${\phi}^1$. \\
We finally obtain for ${\cal{D}}_2$ the only CC: \\
\[ \psi \equiv d_4{\phi}^1-d_3{\phi}^2-d_1{\phi}^3+x^4d_1{\phi}^2-x^3d_1{\phi}^1-{\phi}^1=0    \]

\noindent
{\bf DEFINITION 2.12}: The Janet sequence is said to be {\it locally exact at} $F_r$ if any local section of $F_r$ killed by ${\cal{D}}_{r+1}$ is the image by ${\cal{D}}_r$ of a local section of $F_{r-1}$. It is called {\it locally exact} if it is locally exact at each $F_r$ for $0\leq r \leq n$. The Poincar\'{e} sequence is locally exact, that is a closed form is locally an exact form but counterexamples may exist ([18], p 373).\\

Equivalently, we have the involutive {\it first Spencer operator} $D_1:C_0=R_q\stackrel{j_1}{\rightarrow}J_1(R_q)\rightarrow J_1(R_q)/R_{q+1}\simeq T^*\otimes R_q/\delta (g_{q+1})=C_1$ of order one induced by $D:R_{q+1}\rightarrow T^*\otimes R_q$. Introducing the {\it Spencer bundles} $C_r={\wedge}^rT^*\otimes R_q/{\delta}({\wedge}^{r-1}T^*\otimes g_{q+1})$, the first order involutive ($r+1$)-{\it Spencer operator} $D_{r+1}:C_r\rightarrow C_{r+1}$ is induced by $D:{\wedge}^rT^*\otimes R_{q+1}\rightarrow {\wedge}^{r+1}T^*\otimes R_q:\alpha\otimes {\xi}_{q+1}\rightarrow d\alpha\otimes {\xi}_q+(-1)^r\alpha\wedge D{\xi}_{q+1}$ and we obtain the canonical {\it linear Spencer sequence} ([17], p 150):\\
\[    0 \longrightarrow \Theta \stackrel{j_q}{\longrightarrow} C_0 \stackrel{D_1}{\longrightarrow} C_1 \stackrel{D_2}{\longrightarrow} C_2 \stackrel{D_3}{\longrightarrow} ... \stackrel{D_n}{\longrightarrow} C_n\longrightarrow 0  \]
\noindent
as the canonical Janet sequence for the first order involutive system $R_{q+1}\subset J_1(R_q)$.\\

The canonical Janet sequence and the canonical Spencer sequence can be connected by a commutative diagram where the Spencer sequence is induced by the locally exact central horizontal sequence which is at the same time the Janet sequence for $j_q$ and the Spencer sequence for $J_{q+1}(E)\subset J_1(J_q(E))$ ([17], p 153) but this result will not be used in this paper (See [5],[20],[22],[23] for more details on Cosserat and Maxwell equations, see ([16]-[21]) and in particular ([22],[23]) for applications to engineering and mathematical physics).  \\
 
\noindent
{\bf REMARK 2.13}: We shall revisit Example 2.8 in order to explain the word "{\it canonical} " that has been used in the previous definitions. For this, starting with the {\it inhomogeneous system} $y_{33}=u, y_{13}-y_2=v$, we obtain easily the following {\it inhomogeneous involutive system}  with its corresponding board of multiplicative and non-multiplicative variables:  \\
\[  \left\{  \begin{array}{lclcl}
{\Phi}^4 & \equiv & y_{33}  & = &u  \\
{\Phi}^3 & \equiv & y_{23}  & = & d_1u-d_3v   \\
{\Phi}^2 & \equiv & y_{22} & = & d_{11}u-d_{13}v-d_2v \\
{\Phi}^1 & \equiv & y_{13}-y_2 & = & v
\end{array}  
\right.  \fbox{$\begin{array}{lll}
1 & 2 & 3     \\
1 & 2 & \bullet  \\
1 & 2 & \bullet  \\
1 & \bullet & \bullet 
\end{array}  $  }    \]

\noindent  
Using prolongation with respect to the $4$ non-multiplicative variables involved should bring $4$ first order CC for the right members and we could wait for $4$ {\it third order} CC involving $u$ and $v$.\\
 {\it Surprisingly},  we need the only CC $\Psi \equiv d_{33}v-d_{13}u+d_2u=0$ and obtain the differential sequence:\\
\[   0 \longrightarrow \Theta \longrightarrow \circled{\:1\:} \longrightarrow \circled{\:2\:} \longrightarrow \circled{\:1\:}  \longrightarrow 0  \]
as a single CC has no CC for itself (See ([18],p365) for the effective general procedure). \\
Such a differential sequence is quite different from the canonical Janet sequence: \\
\[  0 \longrightarrow \Theta  \longrightarrow \circled{\:1\:} \stackrel{\cal{D}}{\longrightarrow} \circled{\:4\:} \stackrel{{\cal{D}}_1}{\longrightarrow} \circled{\:4\:} \stackrel{{\cal{D}}_2}{\longrightarrow} \circled{\:1\:} \longrightarrow 0   \]
which is the only sequence that can provide the Spencer sequence as we already said and could not be obtained by simply using Gr\"{o}bner bases. This remark will become essential in mathematical physics (foundations of continuum mechanics, gauge theory, general relativity) where only involutive operators must be used ([20],[22],[23]). We also check that the {\it Euler-Poincar\'e characteristic}, namely the alternate sum of the circled dimensions of the vector bundles involved, does not depend on the differential sequence used as we get $1-2+1=1-4+4-1=0$ (See [18], p 378).\\
\noindent
In the same spirit, using certain parametric jet variables as new unknowns, we may set $z^1=y, z^2=y_1, z^3=y_2, z^4=y_3$ in order to obtain the following involutive first order system with no zero order equation:  \\
\[  \left \{ \begin{array}{lcc}
class \hspace{1mm}  3    & \hspace{1cm}                 &       d_3z^1-z^4=0, d_3z^2-z^3=0, d_3z^3=0,d_3z^4=0  \\
class \hspace{1mm}  2    & \hspace{1cm}                 &       d_2z^1-z^3=0, d_2z^2-d_1z^3=0, d_2z^3=0, d_2z^4=0  \\
class \hspace{1mm}  1    & \hspace{1cm}                 &                      d_1z^1-z^2=0, d_1z^4-z^3=0   
 \end{array} 
\right.   \fbox{  $   \begin{array}{lll}
1  & 2  &  3  \\
1  &  2 &  \bullet  \\
1  &  \bullet  &  \bullet  
\end{array}  $  }   \]

\noindent
where we have separated the classes. Contrary to what could be believed, this operator does not describe the Spencer sequence that could be obtained from the previous Janet sequence. Indeed, introducing the trivial vector bundle $E$ with local coordinates $(x^1,x^2,x^3,y)$, it follows that $J_1(E)$ has local coordinates $(x^1,x^2,x^3,z^1,z^2,z^3, z^4)$. Now, the involutive system $R_2\subset J_2(E)\subset J_1(J_1(E))$ with involutive symbol $g_2\subset S_2T^*\otimes E\subset T^*\otimes T^* \otimes E\subset T^* \otimes J_1(E)$ projects onto $J_1(E)$ but $dim(R_2)=6$ because $par(R_2)=\{y,y_1,y_2,y_3,y_{11},y_{12}\}$ while we have only $4$ unknowns $(z^1,z^2,z^3,z^4)$. Nevertheless, as $R_2$ projects onto $J_1(E)$, we may construct a canonical Janet sequence {\it for this new system} where the successive Janet bundles involved will be the Spencer bundles $C_r={\wedge}^rT^*\otimes J_1(E)/\delta ({\wedge}^{r-1}T^*\otimes g_2)$ {\it with a shift by one step in the numbering of the bundles} as now $C_0=J_1(E)$ and the successive operators are induced by the composition of the inclusion $R_2\subset J_2(E)$ with the Spencer operator $D:J_2(E) \longrightarrow T^*\otimes J_1(E)$ as in ([17],p144,150) or ([18],p356). In any case, it is essential to notice that, both in the canonical Spencer sequence and in the canonical Janet sequence, {\it any intermediate operator can be constructed explicitely without knowing the previous ones}.  \\
 
 \noindent
 {\bf EXAMPLE 2.14}: With $m=1,n=4,q=2$, one could treat similarly the involutive system:  $y_{44}=0, y_{34}=0, y_{33}=0, y_{24}-y_{13}=0$ with one equation of class $4$, two equations of class $3$ and one equation of class $2$.  \\
 
 \noindent
{\bf EXAMPLE 2.15}: Coming back to the involutive system of Example 2.9 with variable coefficients, we let the reader prove that the Janet sequence is:   \\
\[  0 \longrightarrow \Theta \longrightarrow \circled{\:1\:} \stackrel{\cal{D}}{\longrightarrow} \circled{\:3\:} \stackrel{{\cal{D}}_1}{\longrightarrow} \circled{\:3\:} \stackrel{{\cal{D}}_2}{\longrightarrow} \circled{\:1\:} \longrightarrow 0  \] 

\noindent
{\bf EXAMPLE 2.16}: Let us finally consider the following involutive system of PD equations with two independent variables $(x^1,x^2)$ and three unknowns $(y^1,y^2,y^3)$, where again $a$ is an arbitrary constant parameter and we have set for simplicity $y^k_i=d_iy^k$:      \\
\[  \left \{  \begin{array}{ll}
y^2_2+y^2_1+y^3_2-y^3_1-ay^3  &  =  0  \\
y^1_2-y^2_1-y^3_2-y^3_1-ay^3   &  =  0  \\
y^1_1-y^2_1-2y^3_1                       &  =  0
\end{array}  
 \right.  \fbox{ $ \begin{array}{ll}
 1 & 2 \\
 1 & 2 \\
 1 & \bullet 
 \end{array}  $ }  \]

\noindent 
Then the corresponding Janet sequence is:   \\
\[   0 \longrightarrow \Theta \longrightarrow \circled{\:3\:} \stackrel{\cal{D}}{\longrightarrow} \circled{\:3\:} \stackrel{{\cal{D}}_1}{\longrightarrow} \circled{\:1\:} \longrightarrow 0  \]

\noindent
Finally, setting ${\bar{y}}^1=y^1-y^3, {\bar{y}}^2=y^2+y^3$,  we obtain the new first order involutive system:\\
\[  \left\{  \begin{array}{ll}
{\bar{y}}^2_2-{\bar{y}}^2_1-ay^3 & =0 \\  
{\bar{y}}^1_2-{\bar{y}}^2_1-ay^3 & =0 \\
{\bar{y}}^1_1-{\bar{y}}^2_1=0   
\end{array} 
\right.  \fbox{  $ \begin{array}{ll}
1 & 2 \\
1 & 2 \\
1 & \bullet
\end{array}  $ }  \]
with two equations of class 2 and one equation of class 1 in which $y^3$ surprisingly no longer appears.     \\

If ${\chi}_1, ... , {\chi}_n$ are $n$ algebraic indeterminates or, in a more intrinsic way, if $\chi={\chi}_idx^i\in T^*$ is a covector and ${\cal{D}}:E \longrightarrow F:\xi \longrightarrow a^{\tau\mu}_k(x){\partial}_{\mu}{\xi}^k(x)$ is a linear {\it involutive} operator of order $q$, we may introduce the {\it characteristic matrix} $a(x,\chi)=(a^{\tau\mu}_k(x){\chi}_{\mu}\mid \mu=q)$ and the resulting map ${\sigma}_{\chi}({\cal{D}}):E \longrightarrow F$ is called the {\it symbol} of ${\cal{D}}$ at $\chi$. Then there are two possibilities:  \\
\noindent
$\bullet$ If $max_{\chi}rk({\sigma}_{\chi}({\cal{D}})< m \Leftrightarrow {\alpha}^n_q>0$: the characteristic matrix fails to be injective for any covector.\\
\noindent
$\bullet$ If $max_{\chi}rk({\sigma}_{\chi}({\cal{D}})= m\Leftrightarrow {\alpha}^n_q=0$: the characteristic matrix fails to be injective if and only if all the determinants of the $m\times m$ submatrices vanish. However, one can prove that this algebraic ideal $\mathfrak{a} \in K[\chi]$ is not intrinsically defined and must be replaced by its radical $rad(\mathfrak{a})$ made by all polynomials having a power in $\mathfrak{a}$. This radical ideal is called the {\it characteristic ideal} of the operator.\\

\noindent
{\bf DEFINITION 2.17}: For each $x\in X$, the algebraic set defined by the characteristic ideal is called the {\it characteristic set} of $\cal{D}$ at $x$ and $V={\cup}_{x\in X}V_x$ is called the {\it characteristic set} of $\cal{D}$.  \\

One has the following important theorem ([18],[29]) that will play an important part later on:  \\

\noindent
{\bf THEOREM 2.18}: (Hilbert-Serre) The {\it dimension} $d(V)$ of the characteristic set, that is the maximum dimension of the irreducible components, is equal to the number of non-zero characters while the {\it codimension} $cd(V)= n-d(V)$ is equal to the number of zero characters, that is to the number of "{\it full} " classes in the board of multiplicative variables of an involutive system.   \\

\noindent
{\bf EXAMPLE 2.19}: Coming back to Remark 2.12, we obtain $\mathfrak{a}=(({\chi}_3)^2, {\chi}_2{\chi}_3,({\chi}_2)^2, {\chi}_1{\chi}_3) \Longrightarrow rad(\mathfrak{a})=({\chi}_2,{\chi}_3)$ and thus $cd(V)=2$. However, if we take only into account Example 2.8, we should only get the radical ideal $({\chi}_3)$ and the wrong result $cd(V)=1$. The reason for using the radical can be seen from the equivalent first order system that shoul provide $\mathfrak{b}=(({\chi}_3)^4, ... )$ with homogeneous polynomials of degree $4$ and thus $\mathfrak{b}\subset \mathfrak{a}$ with a strict inclusion though $rad(\mathfrak{a})=rad(\mathfrak{b})$. A similar situation can be obtained with Examples 1.1 and 2.9.   \\

\noindent
{\bf 3)  TOOLS FROM MODULE THEORY} :       \\

We may roughly say that, if a reader familiar with Gr\"{o}bner bases ([4],[6]) and computer algebra looks at the previous section, he will feel embarassed because he will believe that "{\it intrinsicness is always competing with complexity} " as can be seen from Examples 2.8 + 2.12. However, even if he admits that it may be useful to have intrinsic and thus canonical procedures, then looking at the existing literature on differential modules ([1],[9],{12]), he will really feel to be on another planet as the main difficulty involved in the theory of differentia modules is to understand why and where formal integrability and involution become essential tools to apply quite before dealing with the homological background of "{\it algebraic analysis} " involving extension modules. This is the main reason for which the case of variable coefficients is rarely treated "{\it by itself} " always refering to Weyl algebras for examples and the main difficulty we found when writing ([18], in particular Chapter IV). The central concept, essential for applications but well hidden in the literature dealing with filtred modules ([14],p 383) and totally absent from the use of Gr\"{o}bner bases because it amounts to formal integrability by duality, is that of a "{\it strict morphism} ". Accordingly, the purpose of this section will be to explain why such a definition, which seems to be purely technical, will be so important for studying extension modules and purity.   \\

If $P=a^{\mu}d_{\mu}\in D=K[d]$, the highest value of ${\mid}\mu {\mid}$ with $a^{\mu}\neq 0$ is called the {\it order} of the {\it operator} $P$ and the ring $D$ with multiplication $(P,Q)\longrightarrow P\circ Q=PQ$ is filtred by the order $q$ of the operators. We have the {\it filtration} $0\subset K=D_0\subset D_1\subset  ... \subset D_q \subset ... \subset D_{\infty}=D$. Moreover, it is clear that $D$, as an algebra, is generated by $K=D_0$ and $T=D_1/D_0$ with $D_1=K\oplus T$ if we identify an element $\xi={\xi}^id_i\in T$ with the vector field $\xi={\xi}^i(x){\partial}_i$ of differential geometry, but with ${\xi}^i\in K$ now. It follows that $D={ }_DD_D$ is a {\it bimodule} over itself, being at the same time a left $D$-module by the composition $P \longrightarrow QP$ and a right $D$-module by the composition $P \longrightarrow PQ$. We define the {\it adjoint} functor $ad:D \longrightarrow D^{op}:P=a^{\mu}d_{\mu} \longrightarrow  ad(P)=(-1)^{\mid \mu \mid}d_{\mu}a^{\mu}$ and we have $ad(ad(P))=P$. It is easy to check that $ad(PQ)=ad(Q)ad(P), \forall P,Q\in D$. Such a definition can also be extended to any matrix of operators by using the transposed matrix of adjoint operators (See [18],[19],[22] for more details and applications to control theory and mathematical physics). \\

Accordingly, if $y=(y^1, ... ,y^m)$ are differential indeterminates, then $D$ acts on $y^k$ by setting $d_iy^k=y^k_i \longrightarrow d_{\mu}y^k=y^k_{\mu}$ with $d_iy^k_{\mu}=y^k_{\mu+1_i}$ and $y^k_0=y^k$. We may therefore use he jet coordinates in a formal way as in the previous section. Therefore, if a system of OD/PD equations is written in the form ${\Phi}^{\tau}\equiv a^{\tau\mu}_ky^k_{\mu}=0$ with coefficients $a\in K$, we may introduce the free differential module $Dy=Dy^1+ ... +Dy^m\simeq D^m$ and consider the differential submodule $I=D\Phi\subset Dy$ which is usually called the {\it module of equations}, both with the {\it differential module} $M=Dy/D\Phi$ or $D$-module and we may set $M={ }_DM$ if we want to specify the ring of differential operators. The work of Macaulay only covers the case $m=1$ with $K$ replaced by $k\subseteq cst(K)$. Again, we may introduce the formal {\it prolongation} with respect to $d_i$ by setting $d_i{\Phi}^{\tau}\equiv a^{\tau\mu}_ky^k_{\mu+1_i}+{\partial}_ia^{\tau\mu}_ky^k_{\mu}$ in order to induce maps $d_i:M \longrightarrow M:{\bar{y} }^k_{\mu} \longrightarrow {\bar{y}}^k_{\mu+1_i}$ by residue if we use to denote the residue $Dy \longrightarrow M: y^k \longrightarrow {\bar{y}}^k$ by a bar as in algebraic geometry. However, for simplicity, we shall not write down the bar when the background will indicate clearly if we are in $Dy$ or in $M$.\\

As a byproduct, the differential modules we shall consider will always be {\it finitely generated} ($k=1,...,m<\infty$) and {\it finitely presented} ($\tau=1, ... ,p<\infty$). Equivalently, introducing the {\it matrix of operators} ${\cal{D}}=(a^{\tau\mu}_kd_{\mu})$ with $m$ columns and $p$ rows, we may introduce the morphism $D^p \stackrel{{\cal{D}}}{\longrightarrow} D^m:(P_{\tau}) \longrightarrow (P_{\tau}{\Phi}^{\tau}):P \longrightarrow P\Phi=P{\cal{D}}$ over $D$ by acting with $D$ {\it on the left of these row vectors} while acting with ${\cal{D}}$ {\it on the right of these row vectors} and the {\it presentation} of $M$ is defined by the exact cokernel sequence $D^p \longrightarrow D^m \longrightarrow M \longrightarrow 0 $. It is essential to notice that the presentation only depends on $K, D$ and $\Phi$ or $ \cal{D}$, that is to say never refers to the concept of (explicit or formal) solutions. It is at this moment that we have to take into account the results of the previous section in order to understant that certain presentations will be much better than others, in particular to establish a link 
with formal integrability and involution. \\

It follows from its definition that $M$ can be endowed with a {\it quotient filtration} obtained from that of $D^m$ which is defined by the order of the jet coordinates $y_q$ in $D_qy$. We have therefore the {\it inductive limit} $0 \subseteq M_0 \subseteq M_1 \subseteq ... \subseteq M_q \subseteq ... \subseteq M_{\infty}=M$ with $d_iM_q\subseteq M_{q+1}$ and $M=DM_q$ for $q\gg 0$ with prolongations $D_rM_q\subseteq M_{q+r}, \forall q,r\geq 0$.  \\

\noindent
{\bf DEFINITION 3.1}: ([14],p 383) If $M$ and $N$ are two differential modules and $f:M \longrightarrow N$ is a morphism over $D$ {\it compatible} with the filtration, that is if $f(M_q)\subset N_q$ with induced morphism $f_q:M_q \longrightarrow N_q$, then $f$ is a {\it strict morphism} if $f_q(M_q)=f(M)\cap N_q, \forall q\geq 0$. \\

Equivalently, chasing in the following diagram: \\
\[   \begin{array}{cccccl}
0 &  & 0 & & 0 &  \\
\downarrow & & \downarrow & & \downarrow &   \\
M_q & \stackrel{f_q}{\longrightarrow} & N_q & \longrightarrow & coker(f_q) & \longrightarrow 0  \\
\downarrow & & \downarrow & & \downarrow &   \\
M & \stackrel{f}{\longrightarrow} & N & \longrightarrow & coker(f) & \longrightarrow 0 
\end{array}  \]
\noindent 
then $f$ is strict if the induced morphism $coker(f_q) \longrightarrow coker(f)$ is a monomorphism $\forall q\geq 0$. \\

\noindent
{\bf DEFINITION 3.2}: An exact sequence of morphisms finishing at $M$ is said to be a {\it resolution} of $M$. If the differential modules involved apart from $
M$ are free, we shall say that we have a {\it free resolution} of $M$. Moreover, a sequence of strict morphisms is called a {\it strict sequence}.  \\

\noindent
{\bf LEMMA 3.3}: If $f$ is a strict morphism as in the last definition, there are exact sequences $0 \longrightarrow coker(f_q) \longrightarrow coker(f_{q+1}) , \forall q\geq 0$.  \\

\noindent
{\it Proof}: As $f$ is compatible with the filtrations and $M_q\subseteq M_{q+1}, N_q\subseteq N_{q+1}$, we have an induced morphism $coker(f_q) \longrightarrow coker(f_{q+1})$. Now, as $f$ is also strict, we have the following commutative and exact diagram:  \\

\[     \begin{array}{rccc}
0\longrightarrow & coker(f_q) & \longrightarrow & coker(f)   \\
  & \downarrow & & \parallel    \\
  0 \longrightarrow & coker(f_{q+1}) & \longrightarrow & coker(f)   
  \end{array}   \]
The lemma finally follows from an elementary chase.  \\
\hspace*{12cm}   Q.E.D.   \\  

Having in mind that $K$ is a left $D$-module with the standard action $(D,K) \longrightarrow K:(d_i,a)\longrightarrow {\partial}_ia$ and that $D$ is a bimodule over itself, {\it we have only two possible constructions}:  \\

\noindent
{\bf DEFINITION 3.4}: We define the {\it system} $R=hom_K(M,K)=M^*$ and set $R_q=hom_K(M_q,K)=M_q^*$ as the {\it system of order} $q$. We have the {\it projective limit} $R=R_{\infty} \longrightarrow ... \longrightarrow R_q \longrightarrow ... \longrightarrow R_1 \longrightarrow R_0$. It follows that $f_q\in R_q:y^k_{\mu} \longrightarrow f^k_{\mu}\in K$ with $a^{\tau\mu}_kf^k_{\mu}=0$ defines a {\it section at order} $q$ and we may set $f_{\infty}=f\in R$ for a {\it section} of $R$. For an arbitrary differential field $K$, {\it such a definition has nothing to do with the concept of a formal power series solution} ({\it care}).\\

\noindent
{\bf DEFINITION 3.5}: We may define the right differential module $hom_D(M,D)$.  \\

\noindent
{\bf PROPOSITION 3.6}: When $M$ is a left $D$-module, then $R$ is also a left $D$-module. \\

\noindent
{\it Proof}: As $D$ is generated by $K$ and $T$ as we already said, let us define:  \\
\[  (af)(m)=af(m), \hspace{4mm} \forall a\in K, \forall m\in M \]
\[ (\xi f)(m)=\xi f(m)-f(\xi m), \hspace{4mm} \forall \xi=a^id_i\in T,\forall m\in M  \]
In the operator sense, it is easy to check that $d_ia=ad_i+{\partial}_ia$ and that $\xi\eta - \eta\xi=[\xi,\eta]$ is the standard bracket of vector fields. We finally 
get $(d_if)^k_{\mu}=(d_if)(y^k_{\mu})={\partial}_if^k_{\mu}-f^k_{\mu +1_i}$ and thus recover {\it exactly} the Spencer operator of the previous section though {\it this is not evident at all}. We also get $(d_id_jf)^k_{\mu}={\partial}_{ij}f^k_{\mu}-{\partial}_if^k_{\mu+1_j}-{\partial}_jf^k_{\mu+1_i}+f^k_{\mu+1_i+1_j} \Longrightarrow d_id_j=d_jd_i, \forall i,j=1,...,n$ and thus $d_iR_{q+1}\subseteq R_q\Longrightarrow d_iR\subset R$ induces a well defined operator $R\longrightarrow T^*\otimes R:f \longrightarrow dx^i\otimes d_if$. This result has been discovered (up to sign) by Macaulay  in 1916 ([11]). For more details on the Spencer operator and its applications, the reader may look at ([22],[23]).  \\
\hspace*{12cm}   Q.E.D.  \\

As $D$ is a bimodule over itself, it follows from this proposition that that $D^*=hom_K(D,K)$ is a left $D$-module. Moreover, using Baer's criterion ([28]), it is known that $D^*$ is an injective $D$-module as there is a canonical isomorphism: \\
\[  M^*=hom_K(M,K)\simeq hom_D(M, D^*)   \]
where both sides are well defined ([2], Prop 11, p 18;)([28], p 37).  \\

\noindent
{\bf DEFINITION 3.7}: With any differential module $M$ we shall associate the {\it graded module} $G=gr(M)$ over the polynomial ring $gr(D)\simeq K[\chi]$ by setting $G={\oplus}^{\infty}_{q=0} G_q$ with $G_q=M_q/M_{q+1}$ and we get $g_q=G_q^*$ where the {\it symbol} $g_q$ is defined by the short exact sequences: \\
\[ 0\longrightarrow M_{q-1}\longrightarrow M_q \longrightarrow G_q \longrightarrow 0  \hspace{4mm}  \Longrightarrow \hspace{4mm}  0 \longrightarrow g_q \longrightarrow R_q \longrightarrow R_{q-1} \longrightarrow 0  \]
We have the short exact sequences $0\longrightarrow D_{q-1} \longrightarrow D_q \longrightarrow S_qT \longrightarrow 0 $ leading to $gr_q(D)\simeq S_qT$ and we may set as usual $T^*=hom_K(T,K)$ in a coherent way with differential geometry. Moreover any compatible morphism $f:M \longrightarrow N$ induces a morphism $gr(f):gr(M) \longrightarrow gr(N)$. \\

\noindent
{\bf EXAMPLE 3.8}: If $K=\mathbb{Q}(x), m=1, n=1$, let us consider the system $y_{xxx}-y_x=0$ for which we may exhibit the basis of sections $\{ f=(1,0,0,0,...),f'=(0,1,0,1,...),f"=(0,0,1,0,...)\}$ as in ([11],\S 59,p 67) or ([21]). We obtain $d_xf=0, d_xf'=-f-f", d_xf"=-f'$ and check that all the sections can be generated by a single one, namely $f"$ which describes the power series of $ch(x)-1$. With now $m=2$, let us consider the module defined by the system $y^1_{xx}=0,y^2_x=0$. Setting $y=y^2-xy^1$, we successively get $y_x=-xy^1_x-y^1, y_{xx}=-2y^1_x, y_{xxx}=0 \Longrightarrow y^1=\frac{x}{2}y_{xx}-y_x,y^2=\frac{x^2}{2}y_{xx}-xy_x+y$ and a differential isomorphism with the module defined by the new system $y_{xxx}=0$. All the sections of the second system are easily seen to be generated by the single section $f=(0,0,1,...)$, a result leading to the only generating section $f^1(x)=\frac{x}{2}, f^1_x=-\frac{1}{2}, f^1_{xx}=0, ..., f^2(x)=\frac{x^2}{2}, f^2_x=0,...$ of the initial system but {\it these sections do not describe solutions} because ${\partial}_xf^1-f^1_x=1\neq 0$ and ${\partial}_xf^2-f^2_x=x\neq 0$. {\it We do not know any reference in computer algebra dealing with sections} (See [21] for more details)  \\

Coming back to the presentation of $M$ under study, we notice that the morphism $\cal{D}$ involved is {\it not} compatible unless we shift the index of the filtration by $ord({\cal{D}})=q$. In that case, we obtain $im({\cal{D}})=I\subset Dy$ and may set $I_{q+r}=I\cap D_{q+r}y$ but we have in general $D^p_r{\cal{D}}\subseteq I_{q+r}$ only, that is the equations of order $q+r$ may not be produced by $r$ prolongations only. We have thus obtained: \\

\noindent
{\bf PROPOSITION 3.9}: The morphism induced by $\cal{D}$ is strict if and only if ${\cal{D}}$ is formally integrable. Accordingly, the module version of both the Janet sequence and the Spencer sequence are strictly exact sequences.  \\

\noindent
{\it Proof}: Using $q+r$ instead of $q$ in Lemma 3.3 and applying $hom_K(\bullet,K)$, we obtain the epimorphisms $R_{q+r+1}\longrightarrow R_{q+r}\longrightarrow 0, \forall r\geq 0$. \\
\hspace*{12cm}  Q.E.D.  \\

The reader will find in ([18], IV,3) more details on the relations existing between $G$ and $M$ which are needed in order to study the non-commutative situation, at least when $K$ is a differential field as such a case is hard enough. We obtain in particular the Hilbert polynomial $dim_K(M_{q+r})=dim_K(R_{q+r})= \frac{\alpha}{d!}r^d+ ...$ where $\alpha$ is called the {\it multiplicity} of $M$ and we use to set $cd_D(M)=cd(M)=n-r, rk_D(M)=rk(M)=\alpha$ if $cd(M)=0$ and $0$ otherwise.  \\

\noindent
{\bf EXAMPLE 3.10}: Coming back to the Example 2.8 of Macaulay, we obtain from Remark 2.13 the free resolution $0\longrightarrow D \longrightarrow D^2 \longrightarrow D \longrightarrow M \longrightarrow 0$ but only the morphism on the left is strict as for the morphism on the right we know that its image is indeed $I_2=\{y_{33},y_{23},y_{22},y_{13}-y_2\} $ and not just $\{y_{33},y_{13}-y_2\}$. However, bringing the system to involution, we get the strict free resolution $0 \longrightarrow D \longrightarrow D^4 \longrightarrow D^4 \longrightarrow D \longrightarrow M \longrightarrow 0$ as the module version of the Janet sequence and we let the reader exhibit the module version of the corresponding Spencer sequence as an exercise.  \\

If $M$ is a differential module over the ring $D=K[d]$ of differential operators and $m\in M$, then the differential submodule $Dm\subset M$ is defined by a system of OD or PD equations for one unknown and we may look for its codimension $cd(Dm)$. A similar comment can be done for any differential submodule $M'\subset M$. Sometimes, a single element $m\in M$, called {\it differentially primitive element}, may generate $M$ if $Dm=M$.\\

\noindent
{\bf EXAMPLE 3.11}: With $K=\mathbb{Q}$, let us consider the differential module$M$ defined by the system $y^1_{xx}-y^1=0, y^2_x=0$. Then, setting $y=y^1-y^2$, we get $y_x=y^1_x, y_{xx}=y^1_{xx}=y^1, y_{xx}-y=y^2$ and thus $y_{xxx}-y_x=0$ provides another way to describe $M$ by means of a single element as in ([18], p435). We have the following commutative and exact diagram:  \\
\[  \begin{array}{rcccccl}
0 \longrightarrow & D^2 & \longrightarrow & D^2 & \longrightarrow & M & \longrightarrow 0  \\
   & \downarrow &  & \downarrow & & \parallel  &   \\
0 \longrightarrow & D & \longrightarrow & D  &\longrightarrow  & M & \longrightarrow 0  \\
  {  }  &  {  } &  {  } & { } & { } & { } & { }  \\
  & (P,Q) & \longrightarrow & (P(d^2-1),Qd) &  &  & \\
  & \downarrow & & \downarrow & & &   \\
    & (Pd+Q) &  \longrightarrow & ((Pd+Q)(d^3-d))  & & & 
 \end{array}  \]  

\noindent
If now we consider the differential module $M$ defined by $y^1_{xx}-ay^1=0,y^2_x=0$ where $a$ is a constant parameter, we cannot find a differentially primitive element when $K=\mathbb{Q}$ if $a=0$ but we can when $K=\mathbb{Q}(x)$ for any value of $a$, as in Example 3.8.  \\

We may check the following definition in a constructive way ([27]):  \\

\noindent
{\bf DEFINITION 3.12}: $t_r(M)=\{m\in M \: {\mid} \: cd(Dm)> r \}$ is the greatest differential submodule of $M$ having codimension $> r$.  \\

\noindent
{\bf PROPOSITION 3.13}: $cd(M)=cd(V)=r \Longleftrightarrow {\alpha}^{n-r}_q\neq 0, {\alpha}^{n-r+1}_q= ... ={\alpha}^n_q=0 \Longleftrightarrow t_r(M)\neq M, t_{r-1}(M)= ... =t_0(M)=t(M)=M$ and this intrinsic result can be most easily checked by using the Spencer form of the system defining $M$.  \\

We are now in a good position for defining and studying purity for differential modules. \\

\noindent
{\bf DEFINITION 3.14}: $M$ is $r$-pure $\Longleftrightarrow t_r(M)=0, t_{r-1}(M)=M \Longleftrightarrow cd(Dm)=r, \forall m\in M$. In particular, $M$ is $0$-pure if $t(M)=0$ and, if $cd(M)=r$ but $M$ is not $r$-pure, we may call $M/t_r(M)$ the pure part of $M$. It follows that $t_{r-1}(M)/t_r(M)$ is equal to zero or is 
$r$-pure (See the picture in [18], p 545). Finally, when $t_{r-1}(M)=t_r(M)$, we shall say that there is a {\it gap} in the purity filtration:   \\
\[   0=t_n(M) \subseteq t_{n-1}(M) \subseteq ... \subseteq t_1(M) \subseteq t_0(M)=t(M) \subseteq M     \]

\noindent
{\bf PROPOSITION 3.15}: $t_r(M)$ does not depend on the presentation or on the filtration of $M$.  \\

\noindent
{\bf EXAMPLE 3.16}: If $K=\mathbb{Q}$ and $M$ is defined by the involutive system $y_{33}=0, y_{23}=0, y_{13}=0$, then $z=y_3$ satifies $d_3z=0, d_2z=0, d_1z=0$ and $cd(Dz)=3$ while $z'=y_2$ only satisfies $d_3z'=0$ and $cd(Dz')=1$. We have the purity filtration 
$0 =t_3(M) \subset t_2(M) =t_1(M) \subset t_0(M)=t(M)=M$ with one gap and two strict inclusions.  \\

We now recall the definition of the {\it extension modules} $ext_D^i(M,D)$ that we shall simply denote by $ext^i(M)$ and the way to use their dimension or codimension. We point out once more that these numbers cannot be obtained without bringing the underlying systems to involution in order to get informations on $M$ from informations on $G$. We divide the procedure into four steps that can be achieved by means of computer algebra ([27]): \\

\noindent
$\bullet$ Construct a {\it free resolution} of $M$, say:  \\
\[   ... \longrightarrow F_i \longrightarrow ... \longrightarrow F_1 \longrightarrow F_0 \longrightarrow M \longrightarrow 0  \]

\noindent
$\bullet$ Suppress $M$ in order to obtain the {\it deleted sequence}:  \\
\[     ... \longrightarrow F_i \longrightarrow ... \longrightarrow F_1 \longrightarrow F_0 \longrightarrow 0  \hspace{11mm}  \]

\noindent
$\bullet$ Apply $hom_D(\bullet,D)$ in order to obtain the {\it dual sequence} heading backwards: \\
\[     ... \longleftarrow hom_D(F_i,D) \longleftarrow ... \longleftarrow hom_D(F_1,D) \longleftarrow hom_D(F_0,D) \longleftarrow 0    \]

\noindent
$\bullet$ Define $ext^i(M)$ to be the cohomology at $hom_D(F_i,D)$ in the dual sequence with $ext^0(M)=hom_D(M,D)$.  \\

The following nested chain of difficult propositions and theorems can be obtained, {\it even in the non-commutative case}, by combining the use of extension modules and {\it bidualizing complexes} in the framework of algebraic analysis. The main difficulty is to obtain first these results for the {\it graded module} $G=gr(M)$ by using techniques from commutative algebra before extending them to the {\it filtred module} $M$ as in ([1],[9],[18],[19]).  \\

\noindent
{\bf THEOREM 3.17}: The extension modules do not depend on the resolution of $M$ used.  \\

\noindent
{\bf PROPOSITION 3.18}: Applying $hom_D(\bullet,D)$ provides right $D$-modules that can be transformed to left $D$-modules by means of the {\it side changing functor} and vice-versa. Namely, if $N_D$ is a right $D$-module, then ${}_DN={\wedge}^nT{\otimes}_KN_D$ is the {\it converted left} $D$-module while, if ${}_DN$ is a left $D$-module, then $N_D={\wedge}^nT^*{\otimes}_K{}_DN$ is the {\it converted right} $D
$-module.\\

\noindent
{\bf PROPOSITION 3.19}: Instead of applying $hom_D(\bullet,D)$ and the side changing functor in the module framework, we may use $ad$ in the operator framework. Namely, to any operator ${\cal{D}}:E \longrightarrow F$ we may associate the formal adjoint $ad({\cal{D}}):{\wedge}^nT^*\otimes F^*\longrightarrow {\wedge}^nT^*\otimes E^*$ with the useful though striking relation $rk_D(ad({\cal{D}}))=rk_D({\cal{D}})$.  \\

\noindent
{\bf PROPOSITION 3.20}: $ext^i(M)$ is a torsion module $\forall 1\leq i \leq n$ but $ext^0(M)=hom_D(M,D)$ may not be a torsion module.  \\

\noindent
{\bf EXAMPLE 3.21}: When $M$ is a torsion module, we have $hom_D(M,D)=0$ (exercise). When $n=3$ and the torsion-free module $M$ is defined by the 
formally surjective $div$ operator, the formal adjoint of $div$ is $-grad$ which defines a torsion module. Also, when $n=1$ as in classical control theory, a controllable system allows to define a torsion-free module $M$ which is free in that case and $hom_D(M,D)$ is thus also a free module. \\

\noindent
{\bf THEOREM 3.22}: \hspace{1cm}  $ext^i(M)=0, \forall i\geq n+1$.  \\

\noindent
{\bf THEOREM 3.23}: \hspace{1cm}  $cd(ext^i(M))\geq i$.  \\

\noindent
{\bf PROPOSITION 3.24}: \hspace{1cm}  $ext^i(M)=0, \forall i< cd(M)$.   \\

\noindent
{\bf THEOREM 3.25}: \hspace{1cm}  $cd(M)\geq r \Leftrightarrow ext^i(M)=0, \forall i<r$.  \\

\noindent
{\bf PROPOSITION 3.26}: \hspace{1cm}  $cd(M)=r \Longrightarrow cd(ext^r(M))=r$ and $ext^r(M)$is $r$-pure.  \\

\noindent
{\bf PROPOSITION 3.27}: \hspace{1cm}  $ext^r(ext^r(M))$ is equal to $0$ or is $r$-pure, $\forall 0\leq r \leq n$.  \\

\noindent
{\bf PROPOSITION 3.28}: If we set $t_{-1}(M)=M$, there are exact sequences $\forall 0\leq r \leq n$:   \\
\[       0 \longrightarrow t_r(M) \longrightarrow t_{r-1}(M)  \longrightarrow ext^r(ext^r(M))  \]

\noindent
{\bf THEOREM 3.29}: If $cd(M)=r$, then $M$ is $r$-pure if and only if there is a monomorphism $0 \longrightarrow M \longrightarrow ext^r(ext^r(M))$ of left differential modules.  \\

\noindent
{\bf THEOREM 3.30}: $M$ is pure $\Longleftrightarrow ext^s(ext^s(M))=0 , \forall s\neq cd(M)$.  \\

The last two theorems are known to characterize purity but it is however evident that they are not very useful in actual practice.\\

\noindent
{\bf THEOREM 3.31}: When $M$ is $r$-pure, the characteristic ideal is thus {\it unmixed}, that is a finite intersection of prime ideals having the same codimension $r$ and the characteristic set is {\it equidimensional}, that is the union of irreducible algebraic varieties having the same codimension $r$.  \\

\noindent
{\bf REMARK 3.32}: For the reader knowing more about commutative algebra, we add a few details about the localization used in the {\it primary decomposition} of a module which are not so well known ([3],[18],[21],[31]). For simplicity, setting $k=cst(K)$, we shall denote by $A=k[\chi]$ the polynomial ring isomorphic to $D=k[d]$ and consider a module $M$ over $A$. We denote as usual by $spec(A)$ the set of {\it proper prime ideals} in $A$, by $max(A)$ the subset of {\it maximal ideals} in $A$ and by $ass(M)=\{\mathfrak{p}\in spec(A) {\mid} \exists 0\neq m\in M, \mathfrak{p}=ann_A(m)\}$ the (finite) set 
$\{ {\mathfrak{p}}_1, ... , {\mathfrak{p}}_t\}$ of {\it associated prime ideals}, while we denote by $\{ {\mathfrak{p}}_1, ...{\mathfrak{p}}_s\}$ the subset of {\it minimum associated prime ideals}. It is well known that $M\neq 0 \Longrightarrow ass(M)\neq \emptyset$. We recall that an ideal $\mathfrak{q}\subset A$ is 
$\mathfrak{p}$-{\it primary} if $ab\in \mathfrak{q}, b\notin \mathfrak{q} \Longrightarrow a\in rad(\mathfrak{q})=\mathfrak{p}\in spec(A)$. We say that a module $Q$ is $\mathfrak{p}$-{\it primary} if $am=0, 0\neq m\in Q \Longrightarrow a\in \mathfrak{p}=rad(\mathfrak{q})\in spec(A)$ when $\mathfrak{q}=ann_A(Q)$ or, equivalently, $ass(Q)=\{\mathfrak{p}\}$. Similarly, we say that a module $P$ is $\mathfrak{p}$-{\it prime} if $am=0, 0\neq m\in P \Longrightarrow a\in \mathfrak{p}\in spec(A)$ when $\mathfrak{p}=ann_A(P)$. It follows that any $\mathfrak{p}$-prime or $\mathfrak{p}$-primary module is $r$-pure with $n-r=trd(A/\mathfrak{p})$, a result generalizing ([11],\S 4, p 43). Accordingly, a module $M$ is $r$-pure if and only if $\mathfrak{a}=ann_A(M)$ admits a primary decomposition $\mathfrak{a}={\mathfrak{q}}_1\cap ... \cap {\mathfrak{q}}_s$ and $rad(\mathfrak{a})={\mathfrak{p}}_1\cap ... \cap {\mathfrak{p}}_s$ with $cd(A/{\mathfrak{p}}_i)=cd(M)=r, \forall i=1, ... ,s$. In that case, the monomorphism $0 \longrightarrow M \longrightarrow {\oplus}_{\mathfrak{p}\in ass(M)}M_{\mathfrak{p}}$ induces a monomorphism $0 \longrightarrow M \longrightarrow Q_1\oplus ... \oplus Q_s$ called {\it primary embedding} where the primary modules $Q_i$ are the images of the localization morphisms $M \longrightarrow M_{{\mathfrak{p}}_i}=S^{-1}M$ with $S=A-\mathfrak{p}$ inducing epimorphisms $M \longrightarrow Q_i \longrightarrow 0$ for $i=1, ... ,s$. Macaulay was only considering the case $M=A/\mathfrak{a}$ with primary decomposition $\mathfrak{a}={\mathfrak{q}}_1\cap ... \cap {\mathfrak{q}}_s$.  \\

\noindent
{\bf EXAMPLE 3.33}: With $k=\mathbb{Q}$ and $n=3$, then $\mathfrak{a}=rad(\mathfrak{a})=({\chi}_1,{\chi}_2{\chi}_3)=({\chi}_1,{\chi}_2)\cap ({\chi}_1,{\chi}_3)$ is unmixed and $M=A/\mathfrak{a}$ is $2$-pure while $\mathfrak{a}=rad(\mathfrak{a})=({\chi}_1{\chi}_2,{\chi}_1{\chi}_3)=({\chi}_1)\cap ({\chi}_2,{\chi}_3)$ is mixed, though an intersection of two minimum prime ideals and $M=A/\mathfrak{a}$ is not $1$-pure. On the contrary, if one has the primary decomposition $\mathfrak{a}=(({\chi}_1)^2, {\chi}_1{\chi}_2, {\chi}_1{\chi}_3, {\chi}_2{\chi}_3)=({\chi}_1, {\chi}_2)\cap ({\chi}_1, {\chi}_3)\cap ({\chi}_1,{\chi}_2, {\chi}_3)^2={\mathfrak{p}}_1\cap {\mathfrak{p}}_2\cap {\mathfrak{m}}^2$ and $M=A/\mathfrak{a}$, then $ass(M)=\{{\mathfrak{p}}_1,{\mathfrak{p}}_2, \mathfrak{m}\}$ with ${\mathfrak{p}}_i \subset \mathfrak{m}$ for $i=1,2$, though $rad(\mathfrak{a})={\mathfrak{p}}_1\cap {\mathfrak{p}}_2$ as before. In that case, there is an embedding $0 \longrightarrow M \longrightarrow Q_1\oplus Q_2 \oplus Q_3$ where $Q_i=A/{\mathfrak{p}}_i$ is the image of the localization morphism $M \longrightarrow M_{{\mathfrak{p}}_i}$ for $i=1,2$ because ${\mathfrak{p}}_1$ is killed by ${\chi}_3\in A-{\mathfrak{p}}_1$, ${\mathfrak{p}}_2$ is killed by 
${\chi}_2\in A-{\mathfrak{p}}_2$ and $Q_3=A/{\mathfrak{m}}^2$ is $\mathfrak{m}$-primary because $rad({\mathfrak{m}}^2)=\mathfrak{m}\in max(A)$. We have also an embedding $0 \longrightarrow M \longrightarrow M_{{\mathfrak{p}}_1}\oplus M_{{\mathfrak{p}}_2}\oplus M_{\mathfrak{m}}$ but no element of the multiplicative set $A-\mathfrak{m}=\{ 1+ a {\mid} a\in \mathfrak{m}\}$ can kill any element of $M$ and the image of $M$ into $M_{\mathfrak{m}}$ is thus isomorphic to $M$ which is not a primay module. It is important to notice that the example of Macaulay $\mathfrak{q}=(({\chi}_3)^2, {\chi}_2{\chi}_3, ({\chi}_2)^2, {\chi}_1{\chi}_3-{\chi}_2)$ provides a $\mathfrak{p}$-primary module $A/{\mathfrak{q}}$ with $\mathfrak{p}=({\chi}_3, {\chi}_2)$ even though the annihilating ideal of $G=gr(M)$ is the homogeneous ideal $\mathfrak{a}=(({\chi}_3)^2, {\chi}_2{\chi}_3, ({\chi}_2)^2, {\chi}_1{\chi}_3)=(({\chi}_2)^2, {\chi}_3)\cap ({\chi}_1,{\chi}_2,{\chi}_3)^2$ which is a mixed ideal because $ass(A/\mathfrak{a})=\{ ({\chi}_2,{\chi}_3),({\chi}_1,{\chi}_2,{\chi}_3)\}$. However, we get $rad(\mathfrak{a})=({\chi}_2,{\chi}_3)$ in a coherent way. \\

\noindent
{\bf PROBLEM }: Is it possible to have a test for checking whether a differential module is pure or not without using the previous results ?   \\

\noindent
{\bf 4)  MOTIVATION} :       \\

As we already said in the introduction and in the previous section, a torsion-free module $M$ is $0$-pure because in that case $t_0(M)=t(M)=0$. Accordingly, $M$ can be embedded into a free module $F$ and the inclusion, which may not be strict when $n>1$, provides a parametrization by means of a finite number of potential-like arbitrary functions in the classical language of elasticity (Airy function) or electromagnetism (EM 4-potential). As it is clear that such a situation is only a very particular case of purity, it remains to wonder what can happen for an $r$-pure module whenever $r\geq 1$. One has the following result ([9], [18], compare to [1], p494):    \\

\noindent
{\bf THEOREM  4.1}: If $M$ is an $r$-pure differential module with $r\geq 1$, there exists a differential module $L$ with $pd(L)\leq r$ and an embedding $M\subseteq L$. \\

\noindent
{\it Proof}: First of all we notice that we have $r>0$ and thus any element $m\in M$ is surely a torsion element because $cd(Dm)>0$, that is $M=t(M)$ is a torsion module with $ext^0(M)=hom_D(M,D)=0$ because $D$ is an integral domain. Let now $ ... \longrightarrow F_r \longrightarrow ... \longrightarrow F_1 \longrightarrow F_0 \longrightarrow N \longrightarrow 0$ be a free resolution of the right differential module $N=ext^r(M)$. According to Proposition 3.26, we have $cd(N)=r>0$ too and $N$ is also a torsion module with $ext^0(N)=hom_D(N,D)=0$. Applying the functor $hom_D(\bullet,D)$ to the previous sequence or , equivalently, constructing the adjoint sequence in the operator framework while using the fact that $ext^i(N)=0, \forall i<r$ according to Theorem 3.25, we obtain the {\it finite long exact sequence with exactly r morphisms} because $N$ is finitely presented and $ext^r(N)\neq 0$:  \\
\[     0 \longrightarrow hom_D(F_0,D) \longrightarrow ... \longrightarrow hom_D(F_{r-1},D) \longrightarrow hom_D(F_r,D) \longrightarrow L \longrightarrow 0 \]

\noindent
where the {\it left} differential module $L$ is the cokernel of the last morphism on the right. As $hom_D(F,D)$ is free whenever $F$ is free because of the bimodule structure of $D={ }_DD_D$, the corresponding deleted complex is:
\[    0 \longrightarrow hom_D(F_0,D) \longrightarrow ... \longrightarrow hom_D(F_{r-1},D) \longrightarrow hom_D(F_r,D)  \longrightarrow 0 \]

\noindent
Applying again $hom_D(\bullet,D)$ and using the {\it reflexivity} of any free module $F$, that is the isomorphism $hom_D(hom_D(F,D),D)\simeq F$, we obtain the dual sequence:  \\
\[  0 \longrightarrow F_r \longrightarrow ... \longrightarrow F_1 \longrightarrow F_0  \longrightarrow 0  \]

\noindent
and a similar procedure may be followed with operators as we shall see in the next illustrating examples ([18],[27]). This sequence is exact everywhere but at $F_r$ and at $F_0$ where its cohomology is just $N$ by definition, that is to say $ext^r(L)=N=ext^r(M)$. Looking for the cohomology at $hom_D(F_r,D)$ in the sequence obtained by duality from the resolution of $N$ with coboundry module $B_r$ and cocycle module $Z_r$, we obtain the following commutative and exact diagram:   \\
\[  \begin{array}{rcccl}
    &  0  &   &   0  &   \\
    & \downarrow & & \downarrow &  \\
0 \longrightarrow & B_r & = & B_r & \longrightarrow 0  \\
   & \downarrow & & \downarrow &  \\
0 \longrightarrow & Z_r & \longrightarrow & hom_D(F_r,D) &    \\
   &  \downarrow & & \downarrow &    \\
0 \longrightarrow &  ext^r(N) & \longrightarrow & L  &   \\
   & \downarrow & & \downarrow &   \\
   & 0  & & 0 & 
 \end{array}   \]

\noindent
Finally, composing the bottom monomorphism with the monomorphism $0 \longrightarrow M \longrightarrow ext ^r(N) $ provided by Theorem 3.29, we get the desired embedding $M \subseteq L$. It must be noticed that such a procedure can be followed equally well in the commutative and non-commutative framework, that is when $K$ is a field of constants or a true differential field.   \\
\hspace*{12cm}      Q.E.D.  \\

\noindent
{\bf EXAMPLE 4.2}: With $K=\mathbb{Q}, m=1, n=4, q=2$, let us study the $2$-pure differential module $M$ defined by the involutive system:   \\
\[  \left\{  \begin{array}{lclll}
{\Phi}^4 & \equiv & y_{44}  & =0  \\
{\Phi}^3 & \equiv &y_{34}  & = 0   \\
{\Phi}^2 & \equiv & y_{33} & = 0  \\
{\Phi}^1& \equiv & y_{24}-y_{13} & = 0
\end{array}  
\right.  \fbox{$\begin{array}{llll}
1 & 2 & 3  & 4   \\
1 & 2 & 3 &\bullet  \\
1 & 2 & 3 &\bullet  \\
1 &2 &\bullet & \bullet 
\end{array}  $  }    \]

\noindent
From the board of multiplicative variables we may construct at once the Janet sequence:  \\
\[  0 \longrightarrow \Theta \longrightarrow \circled{\:1\:} \stackrel{\cal{D}}{\longrightarrow} \circled{\:4\:} \stackrel{{\cal{D}}_1}{\longrightarrow} \circled{\:4\:}
\stackrel{{\cal{D}}_2}{\longrightarrow} \circled{\:1\:} \longrightarrow 0   \]

\noindent
where ${\cal{D}}_1$ is defined by the involutive system:  \\
\[  \left\{  \begin{array}{lclll}
{\Psi}^4 & \equiv & d_4{\Phi}^3-d_3{\Phi}^4  & =0  \\
{\Psi}^3 & \equiv & d_4{\Phi}^2-d3{\Phi}^3 & = 0   \\
{\Psi}^2 & \equiv & d_4{\Phi}^1-d_2{\Phi}^4+d_1{\Phi}^3 & = 0  \\
{\Psi}^1& \equiv & d_3{\Phi}^1-d_2{\Phi}^3+d_1{\Phi}^2 & = 0
\end{array}  
\right.  \fbox{$\begin{array}{llll}
1 & 2 & 3  & 4   \\
1 & 2 & 3 & 4 \\
1 & 2 & 3 & 4 \\
1 & 2 & 3& \bullet 
\end{array}  $  }    \]

\noindent
and ${\cal{D}}_2$ by the (trivially) involutive system:  \\

\[  \left\{  \begin{array}{lclll}
{\Omega} & \equiv & d_4{\Psi}^1-d_3{\Psi}^2+d_2{\Psi}^4-d_1{\Psi}^3  & =0  \\
\end{array}  
\right.  \fbox{$\begin{array}{llll}
1 & 2 & 3  & 4   \\
\end{array}  $  }    \]

\noindent
We have therefore the resolution:  \\
\[ 0 \longrightarrow D \longrightarrow D^4 \longrightarrow D^4 \longrightarrow D \longrightarrow M  \longrightarrow 0  \]

\noindent
leading to $pd(M)\leq 3$ and the deleted complex is:  \\
\[  0 \longrightarrow D \longrightarrow D^4 \longrightarrow D^4 \longrightarrow D \longrightarrow 0  \]

\noindent
Applying $hom_D(\bullet,D)$ to this sequence, we get the sequence:  \\
\[  0  \longleftarrow D \longleftarrow D^4 \longleftarrow D^4 \longleftarrow D \longleftarrow 0   \]

\noindent
which can be described by the following adjoint sequence:  \\
\[    0  \longleftarrow \circled{\:1\:} \stackrel{ad({\cal{D}})}{\longleftarrow} \circled{\:4\:} \stackrel{ad({\cal{D}}_1)}{\longleftarrow} \circled{\:4\:} \stackrel{ad({\cal{D}}_2)}{\longleftarrow} \circled{\:1\:} \longleftarrow 0  \]

\noindent
which is {\it not} a Janet sequence. As $M$ is a torsion module, using now Theorem 3.25 we get $ext^0(M)=0$, $ext^1(M)=0$ and we check that $N=ext^2(M)\neq 0$. For this, dualizing $\Psi$ by $\lambda$ and $\Omega$ by $\theta$, we have to look for the CC of the inhomogeneous system:  \\
\[ \left\{  \begin{array}{lcrl}
{\Psi}^1 & \longrightarrow &-d_4\theta & ={\lambda}^1  \\
{\Psi}^2 & \longrightarrow & d_3\theta & = {\lambda}^2  \\
{\Psi}^3 & \longrightarrow & d_1\theta & = {\lambda}^3  \\
{\Psi}^4 & \longrightarrow & -d_2\theta& = {\lambda}^4
\end{array}
\right.     \]

\noindent
which are not already provided by the system: \\
\[ \left\{  \begin{array}{lcll}
{\Phi}^1 & \longrightarrow & -d_4{\lambda}^2-d_3{\lambda}^1 & =0 \\
{\Phi}^2 & \longrightarrow & -d_4{\lambda}^3-d_1{\lambda}^1 & =0  \\
{\Phi}^3 & \longrightarrow & -(d_4{\lambda}^4-d_2{\lambda}^1)+(d_3{\lambda}^3-d_1{\lambda}^2) & =0  \\
{\Phi}^4 & \longrightarrow & d_3{\lambda}^4+d_2{\lambda}^2 & =0 
\end{array}  
\right.    \]
\noindent
One can check that the torsion module $N$ can be generated by $\{u=d_2{\lambda}^3+d_1{\lambda}^4, v=d_3{\lambda}^3-d_1{\lambda}^2\}$ satisfying the involutive system:  \\
\[   \left\{  \begin{array}{lclll}
{\phi}^4 & \equiv & d_4u-d_1v  & =0  \\
{\phi}^3 & \equiv & d_4v & = 0   \\
{\phi}^2 & \equiv & d_3u-d_2v & = 0  \\
{\phi}^1& \equiv & d_3v & = 0
\end{array}  
\right.  \fbox{$\begin{array}{llll}
1 & 2 & 3  & 4   \\
1 & 2 & 3 & 4 \\
1 & 2 & 3 &\bullet  \\
1 & 2 & 3 & \bullet 
\end{array}  $  }    \]

\noindent
with the two CC:  \\
\[   \left\{  \begin{array}{lclll}
{\psi}^2 & \equiv & d_4{\phi}^2-d_3{\phi}^4+d_2{\phi}^3-d_1{\phi}^1 & = 0  \\
{\psi}^1& \equiv & d_4{\phi}^1-d_3{\phi}^3 & = 0
\end{array}  
\right.  \fbox{$\begin{array}{llll}
1 & 2 & 3  & 4   \\
1 & 2 & 3 & 4 \\
\end{array}  $  }    \]

\noindent
Accordingly, we have the following strict free resolution of $N$:   \\
\[ 0 \longrightarrow D^2 \longrightarrow D^4 \longrightarrow D^2 \longrightarrow N \longrightarrow 0  \]

\noindent
with deleted complex:  \\
\[  0 \longrightarrow D^2 \longrightarrow D^4 \longrightarrow D^2 \longrightarrow 0  \]

\noindent
Applying $hom_D(\bullet,D)$, we get the desired resolution of $L$, namely:  \\
\[ 0  \longleftarrow L \longleftarrow D^2 \longleftarrow D^4 \longleftarrow D^2 \longleftarrow 0  \]

\noindent
Dualizing $\psi$ by $z$, we finally discover that $L$ is defined by the involutive system:  \\
\[  \left\{ \begin{array}{rcll}
-{\phi}^1 & \longrightarrow & d_4z^1-d_1z^2 & =0 \\
-{\phi}^2 & \longrightarrow & d_4z^2      & =0 \\
-{\phi}^3 & \longrightarrow  &  d_3z^1-d_2z^2 & =0 \\
{\phi}^4 & \longrightarrow & d_3z^2   & =0  
\end{array}
\right.  \fbox{ $ \begin{array}{llll}
1 & 2 & 3 & 4  \\
1 & 2 & 3 & 4 \\
1 & 2 & 3 & \bullet \\
1 & 2 & 3 & \bullet 
\end{array} $} \]

\noindent
and is therefore $2$-pure with $pd(L)\leq2$ and a strict inclusion $M\subset L$ defined by $y=z^1$. \\.  \\

\noindent
{\bf REMARK 4.3}: In this example, we discover that, if $L$ were also $r$-pure, we should therefore have an embedding $0 \longrightarrow L \longrightarrow ext^r(ext^r(L))=ext^r(N)$ and thus an isomorphism $ext^r(N)=L$ leading to an isomorphism $Z_r=hom_D(F_r,D)$ and to $F_{r+1}=0$, as can be checked on this example with $r=2$. It has been a challenge for the author during many months to find the following counter-example showing that, {\it sometimes} $L$ {\it may not even be a torsion module}.  \\

\noindent
{\bf EXAMPLE 4.4}: According to the proof of the theorem, $N=ext^r(M)$ does not depend on the resolution of $M$ used while $L$ does indeed depend on the resolution of $N$ used. Coming back to the system studied in Example 2.8 and Remark 2.12 with $r=2$, we may use the shortest finite free resolution of $M$ already presented, namely $0 \longrightarrow D \longrightarrow D^2 \longrightarrow D \longrightarrow M \longrightarrow 0 $. Therefore, taking the adjoint of the only CC found, we may define $N$ by the system:  \\
\[  \left\{   \begin{array}{rcll}
v &\longrightarrow & d_{33}\theta & =0  \\
-u& \longrightarrow & d_{13}\theta +d_2\theta & =0
\end{array}
\right.   \]

\noindent
and obtain the corresponding involutive system:  \\
\[  \left\{  \begin{array}{lcll}
{\phi}^4 & \equiv & d_{33}\theta  & = 0  \\
{\phi}^3 & \equiv & d_{23}\theta  & = 0  \\
{\phi}^2 & \equiv & d_{22}\theta & = 0 \\
{\phi}^1 & \equiv & d_{13}\theta +d_2 \theta & = 0
\end{array}  
\right.  \fbox{$\begin{array}{lll}
1 & 2 & 3     \\
1 & 2 & \bullet  \\
1 & 2 & \bullet  \\
1 & \bullet & \bullet 
\end{array}  $  }    \]

\noindent
We obtain the first order involutive system of CC: \\
\[  \left\{  \begin{array}{rcll}
{\psi}^4 & \equiv & d_3{\phi}^3-d_2{\phi}^4  & = 0  \\
{\psi}^3 & \equiv & d_3{\phi}^2-d_2{\phi}^3 & = 0  \\
{\psi}^2 & \equiv & d_3{\phi}^1-d_1{\phi}^4-{\phi}^3 & = 0  \\
{\psi}^1 & \equiv & d_2{\phi}^1-d_1{\phi}^3-{\phi}^2 & = 0
\end{array}
\right.  \fbox{$\begin{array}{lll}
1 & 2 & 3  \\
1 & 2 & 3   \\
1 & 2 & 3  \\
1 & 2  & \bullet
\end{array} $ }  \]

\noindent
with the only CC : \hspace{1cm} $ \omega \equiv d_3{\psi}^1-d_2{\psi}^2+d_1{\psi}^4+{\psi}^3 =0 $ .  \\
\noindent
We may therefore introduce in reverse order the corresponding adjoint operators of the ones involved in the Janet sequence we have just constructed:  \\
\[  \left\{   \begin{array}{lcrcl}
{\psi}^4 & \longrightarrow & -d_1\lambda & = & {\nu}^4  \\
{\psi}^3 & \longrightarrow &  \lambda & = & {\nu}^3  \\
{\psi}^2 & \longrightarrow  & d_2 \lambda & = & {\nu}^2  \\
{\psi}^1 & \longrightarrow &  -d_3\lambda & = & {\nu}^1
\end{array}   
\right.   \]

\[   \left\{  \begin{array}{lclcll}
{\phi}^4 & \longrightarrow & {\mu}^4 & \equiv & d_2{\nu}^4+d_1{\nu}^2 & =0  \\
{\phi}^3 & \longrightarrow & {\mu}^3 & \equiv & -(d_3{\nu}^4-d_1{\nu}^1)+(d_2{\nu}^3-{\nu}^2) & =0  \\
{\phi}^2 & \longrightarrow & {\mu}^2 & \equiv & -d_3{\nu}^3-{\nu}^1  & =0  \\
{\phi}^1 & \longrightarrow & {\mu}^1 & \equiv & -d_3{\nu}^2-d_2{\nu}^1 & =0
\end{array}
\right.   \]

\noindent
This last operator is defining $L$ but is not involutive. We have the two torsion elements:  \\
\[    \hspace*{1cm }  {\mu}^5 \equiv d_2{\nu}^3-{\nu}^2 , \hspace{1cm} {\mu}^6\equiv d_1{\nu}^3+{\nu}^4  \]

\noindent
which are generating $ext^2(N)$ and are easily seen to satisfy the involutive system:  \\
\[   d_3 {\mu}^6-{\mu}^5 =0, \hspace{1cm} d_3{\mu}^5=0, \hspace{1cm} d_2{\mu}^6-d_1{\mu}^5=0, \hspace{1cm} d_2{\mu}^5=0  \]

\noindent
because $ d_2{\mu}^5 \equiv  d_2 {\mu}^3 + d_3{\mu}^4 + d_1 {\mu}^1 =0$. Finally, using the first equation, we may eliminate ${\mu}^5$ and identify 
${\mu}^6$ with $y$ because we have indeed $d_{33}{\mu}^6=0, d_{13}{\mu}^6-d_2{\mu}^6=0$ in order to obtain the strict inclusion $M\subset L$. Equivalently, we may also eliminate ${\nu}^1$ and ${\nu}^2$ respectively from  ${\mu}^2$ and ${\mu}^3$ in order to obtain:\\
\[ {\mu}^4 \longrightarrow  d_{13}(d_1{\nu}^3+{\nu}^4)-d_2(d_1{\nu}^3+{\nu}^4)=0, \hspace{1cm} {\mu}^1 \longrightarrow d_{33}(d_1{\nu}^3+{\nu}^4)=0  \]

\noindent
but we may notice that $L$ is not $2$-pure and thus a torsion module because ${\nu}^3$ (similarly ${\nu}^4$) is not {\it by itself} a torsion element of $L$. 
Such a situation is well known in control theory with the SISO (single input $u$, single output $y$) system $\dot{y}-\dot{u}=0$ because $u$ (similarly $y$) is not {\it by itself} a torsion element but $z=y-u$ is a torsion element because $\dot{z}=0$ (See the pages 9 and 10 of the introduction in  [17]  for more details on such a comment). \\

\noindent
{\bf PROBLEM   }: Is it possible to find an analogue of the previous theorem or of the case $r=0$, where $L$ should be also $r$-pure with a free resolution having exactly $r$ morphisms ?.  \\

\noindent  
{\bf 5) ABSOLUTE AND RELATIVE LOCALIZATIONS} :         \\

Surprisingly, the positive answer to such a problem has been given by Macaulay in ([11]) for differential modules defined by systems with constant coefficients and only one unknown. Our purpose in this section is to generalize this resul to arbitrary differential modules defined by systems of PD equations with coefficients in a differential field.  \\
 
 Now we hope that, after reading the previous section, the reader is convinced that the use of extension modules is a quite important though striking tool for studying linear multidimensional systems. Of course, as for any new language, it is necessary to apply it on many explicit examples before being familiar with it. However, it is evident that it should be even more important to have a direct approach allowing to exhibit the purity filtration and, in particular, to recognize whether a differential module is pure or not. The purpose of this section is to combine the {\it module approach} with the {\it system approach}, while taking into account the specific properties of the Spencer form in a way rather similar to the use of the Kalman form of a control system when testing controllability, namely to check that an ordinary differential module is $0$-pure. For this, we shall divide the procedure into a few successive constructive steps that will be illustrated on explicit examples.  \\
 
\noindent
$\bullet$ \hspace{2mm} STEP 1: Whenever a system $R_q\subset J_q(E)$ is given, there is no way to obtain informations on the corresponding module without bringing this system to an involutive or at least formally integrable system by means of prolongations and projections as in the Example 2.8 of Macaulay where only the projection $R^{(2)}_2\subset R_2$ of $R_4$ to $R_2$ is involutive. Of course, an homogeneous system with constant coefficients is {\it automatically} formally integrable and one only needs to use a finite number of prolongations in order to obtain an involutive symbol, though it is known that $2$-acyclicity is sufficient to obtain first order generating CC ([17]). However, it is {\it essential} to notice that it is only with an involutive system that we are sure that the CC system is first order both with the following ones in the Janet sequence. \\

\noindent
{\bf EXAMPLE 5.1}: With $K=\mathbb{Q},m=1, n=3,q=2$, the homogeneous second order systems $y_{33}=0, y_{23}-y_{11}=0,y_{22}=0$ or $y_{33}-y_{11}=0, y_{23}=0, y_{22}-y_{11}=0$ both have a $2$-acyclic symbol $g_3$ of dimension $1$ at order $3$ (exercise) and a trivially involutive symbol $g_4=0$ at order $4$, such a result  leading to only one CC of order $2$ with $cd(M)=3$ in both cases. We let the reader treat the system $y_3=0, y_{12}=0$ similarly and conclude (Hint: $({\chi}_3,{\chi}_1{\chi}_2)=({\chi}_3,{\chi}_1)\cap  ({\chi}_3,{\chi}_2)$ is unmixed). It is however not evident that the homogeneous system $y_{11}=0, y_{12}=0, y_{13}=0, y_{23}=0$ of Example 3.33 is involutive.  \\

Finally, according to section 2 and 3, this first step provides the characters ${\alpha}^1_q\geq  ... \geq {\alpha}^n_q\geq 0$ and the smallest non-zero character $\alpha={\alpha}^{n-r}_q\neq 0$ providing $cd(M)=r$, a result leading at once to $t_r(M)\subset M $ with a strict inclusion while $t_{r-1}(M)= ...= t_0(M)=t(M)=M$. Of course, if $\alpha={\alpha}^n_q\neq 0$, then $M$ cannot be a torsion module and $t(M)\subset M$ with a strict inclusion. The following example proves nevertheless that it is much more delicate to study systems with variable coefficients.  \\

\noindent
{\bf EXAMPLE 5.2}: With $K=\mathbb{Q}(x^2) , n=3, m=1, q=1$, let us consider the differential module $M$ defined by the trivially involutive system $y_3-x^2y_1=0$. We have $cd(M)=1$ but we can only say that $cd(Dz)\geq 1, \forall z\in M$. If we set $z=y_2$, proceeding as in Remark  2.12, we get the involutive system:  \\

\[  \left\{  \begin{array}{lcl}
y_3 & = & x^2z_3-(x^2)^2z_1  \\
y_2 & = & z \\
y_1 & = & z_3-x^2 z_1 
\end{array}
\right. \fbox{ $ \begin{array}{lll}
1 & 2 & 3 \\
1 & 2 & \bullet \\
1 & \bullet & \bullet 
\end{array} $ } \]
  
\noindent
The differential submodule $Dz\subset M$ is defined by the {\it second order} involutive system:  \\
\[ \left\{  \begin{array}{ll}
z_{33}-2x^2z_{13}+(x^2)^2z_{11} & =0 \\
z_{23}-x^2 z_{12}-2z_1 & =0
\end{array}
\right. \fbox{ $ \begin{array}{lll}
1 & 2 & 3  \\
1 & 2 & \bullet 
 \end{array} $ }  \]
  
\noindent
and we get $cd(Dz)=1$ exactly. However, even on such a very elementary example, it is not evident that $t_0(M)=t(M)=M$ is $1$-pure. We also understand that the decoupling system for any autonomous element in engineering sciences, like in magnetohydrodynamics, cannot be studied without these new techniques if we want intrinsic results. Finally, if we denote by $I$ the left ideal of $D=Dy$ generated by $y_3-x^2 y_1$, we notice the relation $ann(G)=({\chi}_3-x^2 {\chi}_1)=gr(I)=rad(gr(I))$. However, we have $ ann(gr(Dz))=(({\chi}_3-x^2{\chi}_1)^2, {\chi}_2({\chi}_3-x^2{\chi}_1))$ with radical equal to the prime ideal $({\chi}_3-x^2 {\chi}_1)$ as before. Hence, in this example, the strict inclusion $Dz\subset M$ does not imply $gr((Dz)\subset gr(M)=G$ because otherwise we should get  $ann(G)\subseteq ann(gr(Dz)$ and this is the reason for which only the radical must be considered as it does not depend on the filtration.\\
  
\noindent
$\bullet$ \hspace{2mm} STEP 2: Once we have obtained $cd(M)=r$, in order to check that $M$ is $r$-pure, it remains to prove that $t_r(M)=0$ as we already know that $t_{r-1}(M)=M$. For this, the second step will be to use the specific properties of the Spencer form $R_{q+1}\subset J_1(R_q)$. More generally, it is possible to use any equivalent involutive first order system of the form $R_1\subset J_1(E)$ with no zero order equations, that is with an induced epimorphism $R_1 \longrightarrow E \longrightarrow 0$ and such that the corresponding differential module is isomorphic and thus identified to the initial module as in Remark 2.13 . We have now the characters ${\alpha}^1_1\geq ... \geq {\alpha}^n_1\geq 0$ and the smallest non-zero character is is still $\alpha={\alpha}^{n-r}_1\neq 0$ providing of course the same codimension $cd(M)=r$ as in the first step. Accordingly, the number $r$ of non-zero characters and the number $r$ of full classes is the same as in the previous step. However, it must be noticed that the filtration may be different and the following example explains once more why only the radical of the characteristic ideal must be used.  \\

\noindent
{\bf EXAMPLE 5.3}: $K=\mathbb{Q}, n=2, m=1, q=2$.For the involutive system $y_{22}=0, y_{12}=0$, the characteristic ideal is $\mathfrak{a}=(({\chi}_2)^2, {\chi}_1{\chi}_2)\Longrightarrow rad(\mathfrak{a})=({\chi}_2) \Longrightarrow r=1$. Setting $z^1=y, z^2=y_1, z^3=y_2$, we get the equivalent first order system $d_2z^3=0, d_2z^2=0, d_2z^1-z^3=0, d_1z^1-z^2=0, d_1z^3=0$ and the polynomial ideal generated by the $3\times 3$ minors of the characteristic matrix is  $\mathfrak{a}= (({\chi}_2)^3,({\chi}_2)^2{\chi}_1, {\chi}_2({\chi}_1)^2)$. Hence the characteristic ideal is $rad(\mathfrak{a}) =({\chi}_2) $ and $r=1$ too.  \\

\noindent
{\bf EXAMPLE 5.4}: For Example 2.8 we may set $z^1=y, z^2=y_1, z^3=y_2, z^4=y_3$ and obtain the first order involutive system :ÊÊ\\
\[ \left\{   \begin{array}{c}
d_3z^4=0, d_3z^3=0, d_3z^2-z^3=0, d_3z^1-z^4  =0  \\
d_2z^4=0, d_2z^3=0, d_2z^2-d_1z^3=0, d_2z^1-z^3=0   \\
d_1z^4-z^3=0, d_1z^1-z^2  =0
\end{array}
\right. \fbox{ $ \begin{array}{lll}
1 & 2 & 3 \\
1 & 2 & \bullet \\
1 & \bullet & \bullet
\end{array} $ }      \]

\noindent
with no zero-order equation. We have ${\alpha}^3_1=4-4=0, {\alpha}^2_1=4-4=0, {\alpha}^1_1=4-2=2  \Longrightarrow r=2$ too. We let the reader treat Example 4.2 similarly and obtain ${\alpha}^4_1=0, {\alpha}^3_1=0, {\alpha}^2_1=2  \Longrightarrow r=2$.  \\

It is at this moment that we discover that such systems have particular properties not held by other systems, apart from the fact that a canonical sequence may be constructed exactly like the Spencer sequence or the first order part of the Janet sequence. \\
  
Shrinking the board of multiplicative variables, we obtain from the definition of involutiveness:  \\  
  
\noindent
{\bf PROPOSITION 5.5}: For an involutive first order system with no zero order equations and solved with respect to the principal ($pri$) first order jets expressed by means of the parametric ($par$) other first order jets, the system obtained by looking only at the PD equations of class $1$+ ... + class $i$ only contains $d_1, ..., d_i$ and is still involutive $\forall i=1, ..., n$, after adopting the ordering $d_{i+1}, ..., d_n, d_1,..., d_i$. \\ 
  
\noindent
{\bf EXAMPLE 5.6}: Looking at Example 2.9, we notice that the systems:  \\
\[  \left\{  \begin{array}{ll}
y_3-x^4 y_1 & =0  \\
y_2 - y_1 & =0 
\end{array}
\right.  \fbox{$\begin{array}{llll}
4 & 1 & 2 & 3  \\
4 & 1 & 2 & \bullet
\end{array} $ }  \]
 
\[   \left\{  \begin{array}{ll}
y_2 - y_1 & =0  
\end{array}
\right.   \fbox{ $ \begin{array}{llll}
3 & 4 & 1 & 2 
\end{array} $ }  \]

\noindent
are both involutive. Also, looking at Example 5.4, we notice that the system:
\[  \left\{  \begin{array}{c}
d_2z^4=0,d_2z^3=0, d_2z^2-d_1z^3=0, d_2z^1 - z^3=0  \\
d_1z^4-z^3=0, d_1z^1-z^2=0
\end{array}
\right.  \fbox{$ \begin{array}{lll}
3 & 1 & 2 \\
3 & 1 & \bullet
\end{array}  $  }   \]

\noindent
is still involutive and we let the reader treat Example 4.2 similarly.  \\

We shall denote the corresponding differential module by $M_{n-i}$ and we notice that $M=M_{0}$ is defined by {\it more equations} than $M_{n-i}$. Accordingly, we have an epimorphism ({\it specialization}) $M_{n-i} \longrightarrow M  \longrightarrow 0  $ and similarly epimorphisms $M_{n-i} \longrightarrow M_{n-i-1} \longrightarrow 0  $. Finally, as $cd(M)=r$, we notice that the classes $n-r+1, ..., n$ are full and we find therefore $({\chi}_n)^m, ... , ({\chi}_{n-r+1})^m$ among the $m\times m$ minors with lower powers of ${\chi}_1, ... , {\chi}_{n-r}$ for the other minors because the numbers of equations of the lower classes are decreasing and thus strictly smaller than $m$. The characteristic ideal is thus $({\chi}_n, ..., {\chi}_{n-r+1})$ if we set ${\chi}_1, ... , {\chi}_{n-r}$ to zero, in a coherent way. Finally, choosing $i=n-r$, we get an epimorphism $M_{r} \longrightarrow M \longrightarrow 0$. The background will always indicate clearly the meaning of the lower index and cannot be confused with the filtration index of $M$.\\
 
\noindent
$\bullet$ \hspace{2mm} STEP 3: We are now in position to study $M_r$ with more details as a system in $n-r$ variables ([11], \S 77, p 86). Its defining system has ${\beta}^{n-r}_1=\beta=m-\alpha$ equations of strict class $n-r$, a smaller number of equations of class $n-r-1$, ... , and eventually an even smaller number of equations of class $1$. Studying this system {\it for itself}, we may look for $t(M_r)$ exactly following the known techniques working for any differential module, in particular double duality as described in section 4. However, if one could find any (relative) torsion element $z\in M_r$, we could project it to an element $z\in M$ and we have $N=Dz\subseteq M$ where we do not put a residue bar on the new $z$ for simplicity. Introducing the respective annihilator ideal $\mathfrak{a}$ and $\mathfrak{b}$ of $M$ and $N$, we should have $\mathfrak{a}\subseteq rad(\mathfrak{a})\subseteq rad(\mathfrak{b})$ as it is the only result not depending on the filtration of the modules. However, we know that $z$ must be killed by at least one operator involving only $d_1, ... , d_{n-r}$, in addition to the operators involving {\it separately} $(d_n)^m+ ... , ... , (d_{n-r-1})^m+ ... $ and we should obtain $t_{r-1}(M)=M$ but $t_r(M)\neq 0$, that is $M$ should not be pure. Hence $M$ is $r$-pure if and only if $M_r$ is torsion-free as a differential module over $K[d_1,...,d_{n-r}]$. In such a case, the system defining $M_r$ can be  parametrized by $\alpha$ arbitrary unknowns among $\{y^1, ... , y^m\}$ by using a so-called {\it minimum parametrization} in the sense of ([24]). In actual practice, as shown on all the examples, one can use a {\it relative localization} with respect to $(d_1, ... , d_{n-r})$ {\it only} by keeping $(d_{n-r+1}, ... , d_n)$ untouched and replacing $(d_1, ... , d_{n-r})$ by $({\chi}_1, ... , {\chi}_{n-r})$ considered as (constant) "parameters" in the language of Macaulay ([11], \S 43, p 45 and \S 77, p 86 with $r$ instead of $n-r$ and a different ordering). Such a method provides therefore a quite useful and simple test for checking purity by linking it to involutivity. {\it An important intermediate result is provided by the next proposition}.  \\

\noindent
{\bf PROPOSITION 5.7}: The partial localization "{\it kills}" the equations of class $1$ up to class $n-r-1$ ({\it care}) and finally only depends on the equations of {\it strict} class $n-r$.  \\

\noindent
{\it Proof}: Instead of using the column $n-r$ in the multiplicative board, we provide the proof when $r=0$ by using the column $n$. Working out as usual the {\it first order} CC, we only look at the $p<m$ CC of class $n$ for the equations ${\Phi}^1, ... , {\Phi}^p$ of class $1$ up to class $n-1$ if we order the $\Phi$ starting from the lower class involved. These $p$ CC will be of a very specific form with a square $p\times p$ operator matrix for $({\Phi}^1, ... , {\Phi}^p) $ with diagonal operators of the form $d_n+ ... $ where the dots denote operators involving only $d_1, ... , d_{n-1}$ in a quasi linear way with coefficients in $K$, the remaining of the matrix only depending on $d_1, ... , d_{n-1}$ for the $({\Phi}^{p+1}, ... , {\Phi}^m)$ of {\it strict class n}. Therefore, if we have $K=k=cst(K)$, the 
{\it absolute localization} is simply done by setting $d_i\longrightarrow {\chi}_i$ and the determinant of the square matrix is equal to $({\chi}_n)^p+ ... $ where the dots denote a polynomial of degree $<p$ in ${\chi}_p$ with coefficients involving only ${\chi}_1, ... , {\chi}_{n-1}$. It follows that each ${\Phi}^1, ... , {\Phi}^p$ can be expressed as a linear combination over $k({\chi}_1, ... , {\chi}_n)$ of the $\Phi$ of strict class $n$. The result is similar for the variable coefficient 
case by using the graded machinery but needs much more work. In any case, setting the $\Phi$ of strict class $n$ equal to zero, we should obtain a zero graded module for the ${\Phi}^1, ... , {\Phi}^p$ which must be eequal to zero too. It must finally be noticed that the first order CC used are in reduced Spencer form as the $d_n$ of the $\Phi$ of strict class $n$ do not appear in the CC we have used and {\it these} $\Phi$ {\it do not appear in the other} CC {\it too}.  \\
\hspace*{12cm}  Q.E.D   \\

\noindent
{\bf EXAMPLE 5.8}: Coming back to Example 2.11, we notice that the only CC is an identity in $(u,v,w)$ and we may forget about $u$ in order to obtain the new system for $(v,w)$:  \\
\[   \left\{  \begin{array}{ll}
{\phi}^3\equiv & d_4v-x^3d_1v -v =0 \\
{\phi}^2\equiv  & d_4w-x^3d_1w-w=0  \\
{\phi}^1\equiv & d_3w-x^4d_1w-d_2v+d_1v=0
\end{array}
\right.  \fbox{  $ \begin{array}{cccc}
1 & 2 & 3 & 4 \\
1 & 2 & 3 & 4 \\
1 & 2 & 3 & \bullet 
\end{array}  $ }   \]

\noindent
defining a new module $M$ with $cd(M)=1$ as the class $4$ is {\it now} full. The same CC as before can be written again in the form:\\
\[ (d_4-x^3d_1-1){\phi}^1=(d_3-x^4d_1){\phi}^2-(d_2-d_1){\phi}^3     \]
that we can localize in $(d_1,d_2,d_3)$ with:  \\
\[  {\phi}^1=0  \Longrightarrow (d_3-x^4d_1)w=(d_2-d_1)v \Longrightarrow v=(d_3-x^4d_1)y, w=(d_2-d_1)y  \]
 \noindent
 We let the reader check that we have indeed:  \\
 \[  \left\{  \begin{array}{rrcrl}
 {\phi}^2\equiv & (d_4-x^3d_1-1)(d_3-x^4d_1)y & = & (d_3-x^4d_1)u =0  & \Longrightarrow  u=0  \\
 {\phi}^3\equiv & (d_4-x^3d_1-1)(d_2-d_1)y & = & (d_2-d_1)u=0 &\Longrightarrow u=0
 \end{array}   \right.  \]
 We finally obtain for $M$ a relative parametrization with the only constraint $u\equiv (d_4-x^3d_1-1)y=0$ in a coherent way with Example 2.11.  \\
 
\noindent
$\bullet$ \hspace{2mm} STEP 4: In this section we come back to the commutative situation with a field $k$ of constants and generalize the results of Macaulay in the following theorem which is recapitulating the results so far obtained.   \\

\noindent
{\bf THEOREM 5.9}: One has the commutative and exact diagram: \\

\[  \begin{array}{rccccrcl}
0\longrightarrow & t(M_r) & \longrightarrow & M_r & \longrightarrow & k({\chi}_1, ... , {\chi}_{n-r}) & \otimes & M_r  \\
     & \downarrow &  & \downarrow &  &   & \downarrow &   \\
     0 \longrightarrow & t_r(M) & \longrightarrow & M & \longrightarrow & k({\chi}_1, ... , {\chi}_{n-r}) & \otimes & M  \\
         & \downarrow &  & \downarrow & & & \downarrow &\\
          & 0  & & 0 & & & 0 & 
          \end{array}  \]
 
\noindent
{\it Proof}: For simplifying the notations in this diagram of modules over $D$, we have identified $M_r$ as a module over $k[d_1, ... , d_{n-r}]$ with $M_r$ as a module over $k[d_{n-r+1}, ... , d_n, d_1, ... , d_{n-r}]=k[d_1, ... ,d_n]=D$ while the localization of $M$ just tells that the coefficients are now in the field $k({\chi}_1, ... ,{\chi}_{n-r})$, exactly following Macaulay. Moreover the central column is exact according to the definitionof $M_r$ and the right column is exact because localization preserves the exactness of a sequence.  \\
For exampe, with $k=\mathbb{Q}, n=2,m=2, q=2, r=1$, the differential module $M$ defined by the involutive system $y_{22}=0, y_{12}=0$ may also be defined by the first order involutive system $z^1=y, z^2=y_1, z^3=y_2 \Longrightarrow d_2z^3=0, d_2z^2=0, d_2z^1-z^3=0, d_1z^3=0, d_1z^1-z^2=0$. Then $M_1$ is defined by the first order system $d_1z^3=0, d_1z^1-z^2=0$ with torsion module $t(M_1)$ generated by $z^3$ and the tensor product of $M$ by 
$k({\chi}_1)$ is defined by $ d_2z^3=0, d_2z^2=0, d_2z^1=0, d_1z^3=0, z^3=0, {\chi}_1z^1-z^2=0$ after division by ${\chi}_1$ in a coherent way with Example 1.1. \\
\hspace*{12cm}  Q.E.D.  \\

Finally, when $M_r$ is torsion-free as a differential module over $k[d_1, ... ,d_{n-r}]$, then $t_r(M)=0$ and we get the following generalization of a result provided by Macalualy ([11], \S 41, p 43):  \\

\noindent
{\bf COROLLARY 5.10}: The differential module $M$ is $r$-pure if and only if  $cd(M)=r$ and there is a monomorphism $0 \longrightarrow M \longrightarrow k({\chi}_1, ... ,{\chi}_{n-r}) \otimes M  $. \\
  
\noindent
$\bullet $\hspace{2mm} STEP 5: The final idea is to embed $M_r$ into a free module over $K[d_1, ..., d_{n-r}]$ in order to parametrize the corresponding system and substitute into the equations of class $n-r+1, ... , n$. However, if we look at Example 1.2, we should find after the substitution ${\Phi}^1\equiv z_{13}=0, {\Phi}^2\equiv z_{23}=0$ with one CC $d_2{\Phi}^1-d_1{\Phi}^2=0$, that is on one side a module $L$ which is not $1$-pure and, on the other side a module $L$ having a finite free resolution with $2$ operators. However, we forgot that $M$, being pure, may be identified with its embedding into its localization. Hence, we get in fact ${\chi}_1z_3=0, {\chi}_2z_3=0$ and thus only $z_3=0$ is providing a convenient {\it parametrizing module} $L$. \\

Our purpose is to explain and illustrate this procedure for finding such an $L$ in the general situation. Again, the main idea will be provided by this example. Indeed, we obtain the only CC $\Psi \equiv d_3{\Phi}^3-d_2{\Phi}^1+d_1{\Phi}^2=0$. Substituting the parametrization, we get of course ${\Phi}^3=0\Longleftrightarrow {\chi}_1y^2={\chi}_2y^1$, that is, among the two unknowns $y^1,y^2$ we are left with only one, say $y^1$ and, similarly, among the two equations ${\Phi}^1,{\Phi}^2$ we are left with only one, say ${\Phi}^1$, because ${\chi}_1{\Phi}^2={\chi}_2{\Phi}^1$ from the CC which is of course compatible with the localization and we choose $z_3=0$ as ${\chi}_1z_3=0 \Longrightarrow z_3=0$.  \\

The general situation may be treated similarly. Indeed, according to the previous step, we are only concerned with the equations of class $n-r+1$, ... , class 
$n$ while the localization has only to do with the $\beta$ equations of strict class $n-r$ ({\it care}) allowing to express $\beta$ unknowns as linear combinations of the $\alpha$ remaining unknowns with coefficients in $k({\chi}_1, ... , {\chi}_{n-r})$. To each such equation are associated exactly $r$ dots and each dot of index $n-r+i$ provides a reduction of the respective equations of class $n-r+i$ for $i=1, ... , r$. It follows that we are left with $\alpha$ equations of each such class. When we "delocalize", replacing ${\chi}_i$ by $d_i$, we have to take into account the need to take out the denominators and may find a few "simplifications" as in the example just considered.. Finally, the maximum number $r-1$ ({\it care again}) of dots found for one equation is obtained for the equations of strict class $n-r+1=n-(r-1)$ and we have thus exhibited a system defining a module $L$ which is $r$-pure and admits a free resolution with exactly $(r-1)+1=r$ operators. In any case, the reader must not forget that the localization of a module is useful only if we already know that this module is torsion-free by means of the double-duality formula $t(M)=ext^1(N)$ given in the introduction.   \\
 
\noindent
{\bf EXAMPLE 5.11}: Let $M$ be defined by the involutive system: \\
\[   \left\{    \begin{array}{ll}
d_3y^4+d_1y^2-d_1y^1 & =0  \\
d_3y^3-d_2y^4+d_1y^2-d_1y^1 & =0  \\
d_3y^2+d_1y^2 & =0 \\
d_3y^1-d_1y^4+d_1y^2 & =0  \\
d_2y^2-d_1y^4+d_1y^1 & =0  \\
d_2y^1-d_1y^3+d_1y^1 & =0
\end{array}  
\right.  \fbox{ $ \begin{array}{ccc}
1 & 2 & 3 \\
1 & 2 & 3 \\
1 & 2 & 3 \\
1 & 2 & 3 \\
1 & 2 & \bullet \\
1 & 2 & \bullet
\end{array}  $ }   \]
with characters ${\alpha}^3_1=0, {\alpha}^2_1=\alpha=2, {\alpha}^1_1=4$. It follows that $cd(M)=1$ as only the class $4$ is full and we obtain the 
following relative localization showing that $M$ is $1$-pure: \\
\[ y^1={\chi}_1y, y^2={\chi}_1z, y^3=({\chi}_1+{\chi}_2)y, y^4={\chi}_1y+{\chi}_2z  \]
Substituting in the four equations of class $3$, we only obtain the two equations: \\
\[   \left\{  \begin{array}{ll}
d_3z+{\chi}_1z & =0 \\
d_3y+({\chi}_1-{\chi}_2)z-{\chi}_1y & =0
\end{array}
\right.   \]
after a division by ${\chi}_1,{\chi}_2$ and ${\chi}_1+{\chi}_2$. The parametrizing module $L$ is thus defined by the two equations: \\
\[  \left\{   \begin{array}{ll}
d_3z+d_1z & =0  \\
d_3y+(d_1-d_2)z-d_1y & =0
\end{array} 
\right.    \]
which are differentially independent and we have the relative parametrization:  \\ 
\[  y^1=d_1y,\hspace{3mm} y^2=d_1z,\hspace{3mm} y^3=(d_1+d_2)y,\hspace{3mm} y^4=d_1y+d_2z  \]
Finally, $M\subset L\Longrightarrow ass(M)\subset ass(L)\Longrightarrow ass(M)=\{(d_3+d_1),(d_3-d_1)\} \Longrightarrow ann_D(M)=(d_3+d_1)\cap (d_3-d_1)$, a striking result showing that $M$ can be embedded into the direct sum of two primary differential modules according to Remark 3.32 
(See [18] for more details).   \\

\noindent
{\bf EXAMPLE 5.12}: With $k=\mathbb{Q}, m=1, n=3$, let us consider the polynomial map ${\chi}_1=u^5, {\chi}_2=u^3, {\chi}_3=u^4$ as in ([11], p 53). We have the exact sequence $0\longrightarrow \mathfrak{p} \longrightarrow k[\chi] \longrightarrow k[u]\subset k(u)$ showing that $\mathfrak{p}=
(({\chi}_2)^2({\chi}_3)-({\chi}_1)^2, ({\chi}_2)^3-{\chi}_1{\chi}_3, ({\chi}_3)^2-{\chi}_1{\chi}_2)$ is a prime ideal ([18], p 126). The corresponding prime differential module $M$ is defined by the involutive system: \\
\[  \left\{  \begin{array}{ll}
y_{333}-y_{123} & =0 \\
y_{233}-y_{122} & =0 \\
y_{223}-y_{11} & =0 \\
y_{222}-y_{13} & =0 \\
y_{133}-y_{112} & =0 \\
y_{33}-y_{12} & =0 
\end{array}
\right. \fbox{ $ \begin{array}{lll}
1 & 2 & 3 \\
1 & 2 & \bullet \\
1 & 2 & \bullet \\
1 & 2 & \bullet \\
1 & \bullet & \bullet \\
\bullet & \bullet & \bullet 
\end{array}  $ } \]
and is $2$-pure. The localized system is finite type over $k({\chi}_1)[d_2, d_3]$ with $par=\{y,y_2,y_3,y_{22},y_{23}\}$. One can prove that ${\mathfrak{p}}^2$ is not a primary ideal even though $rad({\mathfrak{p}}^2)=\mathfrak{p}$. \\
 
\noindent
{\bf EXAMPLE 5.13}: Similarly but now with $k=\mathbb{Q}, m=1, n=4$, let us consider the polynomial map ${\chi}_1=uv, {\chi}_2=u, {\chi}_3=uv^3, {\chi}_4=uv^4$ as in ([11], p 47). We have the exact sequence $0 \longrightarrow \mathfrak{p} \longrightarrow k[\chi] \longrightarrow k[u,v]\subset k(u,v) $ showing that $\mathfrak{p}=({\chi}_2{\chi}_4-{\chi}_1{\chi}_3, ({\chi}_1)^3-({\chi}_2)^2{\chi}_3, ({\chi}_3)^3-{\chi}_1({\chi}_4)^2, ({\chi}_1)^2{\chi}_4-{\chi}_2({\chi}_3)^2)$ is a prime ideal. It is not evident at all that the corresponding prime differential module $M$ can be defined by the homogeneous involutive system (exercise): \\
\[   \left\{  \begin{array}{ll}
y_{444}-y_{224}-y_{134}-y_{123} &=0 \\
y_{344}-y_{111} & =0 \\
y_{334}-y_{114}-y_{112} & =0 \\
y_{333}-y_{124}-y_{122}-y_{113} & =0 \\
y_{244}+y_{224}-y_{123} & =0 \\
y_{234}-y_{133}+y_{111} & =0  \\
y_{144}+y_{124}-y_{113} & =0 \\
y_{44}+y_{24}-y_{13} & =0
\end{array}
\right. \fbox{ $ \begin{array}{llll}
1 & 2 & 3 & 4 \\
1 & 2 & 3 & \bullet \\
1 & 2 & 3 & \bullet \\
1 & 2 & 3 & \bullet \\
1 & 2 & \bullet & \bullet \\
1 & 2 & \bullet & \bullet \\
1 & \bullet & \bullet & \bullet \\
\bullet & \bullet & \bullet & \bullet 
\end{array}  $ }  \]
and is thus also $2$-pure. The localized system is finite type over $k({\chi}_1,{\chi}_2)[d_3,d_4]$ with $par=\{ y, y_3, y_4, y_{33}\}$ and we have for example ${\chi}_2y_{34}-{\chi}_1y_{33}+({\chi}_1)^3y=0$ in a coherent way with the comments of Macaulay in ([11], \S 78, p 88, formula (A) and \S 88,89, p 98). \\

\noindent
{\bf 6)  CONCLUSION} :           \\

In 1916, F.S. Macaulay discovered a new localization technique for studying {\it unmixed polynomial ideals}. We have been able to generalize this procedure for studying {\it pure differential modules}, obtaining in particular a kind of {\it relative parametrization} generalizing the {\it absolute parametrization} already known for torsion-free modules and equivalent to controllability in classical control theory. In the language of multidimensional systems theory, which is more intuitive, instead of using arbitrary potential-like functions for the parametrization, the idea is now to use potential-like functions which must satisfy a kind of {\it minimum differential constraint} limiting, in some sense, the number of independent variables appearing in these functions, in a way similar to the situation met in the Cartan-Kh\"{a}ler theorem of analysis. For such a purpose, we have exhibited new links between purity and involutivity, providing also a new insight into the primary decomposition of modules and ideals by means of tools from the formal theory of linear multidimensional systems.  \\

    \noindent
{\bf 7)  BIBLIOGRAPHY} :     \\

\noindent
[1] J.E. BJORK: Analytic D-modules and Applications, Kluwer, 1993.\\
\noindent
[2] N. BOURBAKI: Alg\`{e}bre homologique, chap. X, Masson, Paris, 1980.  \\
\noindent
[3] N. BOURBAKI: Alg\`{e}bre commutative, chap. I-IV, Masson, Paris, 1985.  \\
\noindent
[4] B. BUCHBERGER: Gr\"{o}bner bases: an Algorithmic Methods in Polynomial Ideal Theory, in: Recent Trends in Multidimensional System Theory, 
N.K. Bose ED., Reidel, Dordrecht, 1985, 184-232.\\ 
\noindent
[5] E. COSSERAT, F. COSSERAT: Th\'{e}orie des Corps D\'{e}formables, Hermann, Paris, 1909.\\
\noindent
[6] W. GR\"{O}BNER: \"{U}ber die Algebraischen Eigenschaften der Integrale von Linearen Differentialgleichungen mit Konstanten Koeffizienten, 
Monatsh. der Math., 47, 1939, 247-284.\\
\noindent
[7] M. JANET: Sur les Syst\`emes aux d\'eriv\'ees partielles, Journal de Math., 8, 3, 1920, 65-151.\\
\noindent
[8] E.R. KALMAN, Y.C. YO, K.S. NARENDA: Controllability of Linear Dynamical Systems, Contrib. Diff. Equations, 1, 2, 1963, 189-213.\\
\noindent
[9] M. KASHIWARA: Algebraic Study of Systems of Partial Differential Equations, M\'emoires de la Soci\'et\'e Math\'ematique de France 63, 1995, 
(Transl. from Japanese of his 1970 Master's Thesis).\\
\noindent
[10] E. KUNZ: Introduction to Commutative Algebra and Algebraic Geometry, Birkh\"{a}user, 1985.\\
\noindent
[11] F.S. MACAULAY: The Algebraic Theory of Modular Systems, Cambridge Tracts, vol. 19, Cambridge University Press, London, 1916. Stechert-Hafner Service Agency, New-York, 1964.\\
\noindent
[12] P. MAISONOBE, C. SABBAH: D-Modules Coh\'erents et Holonomes, Travaux en Cours, 45, Hermann, Paris, 1993.\\
\noindent
[13] D.G. NORTHCOTT: An Introduction to Homological Algebra, Cambridge University Press, Cambridge, 1966.   \\
\noindent
[14] D.G. NORTHCOTT: Lessons on Rings, Modules and Multiplicities, Cambridge University Press, Cambridge, 1968.   \\
\noindent
[15] V.P. PALAMODOV: Linear Differential Operators with Constant Coefficients, Grundlehren der Mathematischen Wissenschaften 168, Springer, 1970.\\
\noindent
[16] J.-F. POMMARET,:Dualit\'e Diff\'erentielle et Applications, C. R. Acad. Sci. Paris, 320, S\'erie I, 1995, 1225-1230.\\
\noindent
[17] J.-F. POMMARET: Partial Differential Equations and Group Theory: New Perspectives for Applications, Kluwer, 1994.\\
\noindent 
[18] J.-F. POMMARET: Partial Differential Control Theory, Kluwer, 2001.\\
\noindent
[19] J.-F. POMMARET: Algebraic Analysis of Control Systems Defined by Partial Differential Equations, in Advanced Topics in Control Systems Theory, Lecture Notes in Control and Information Sciences LNCIS 311, Chapter 5, Springer, 2005, 155-223.\\
\noindent
[20] J.-F. POMMARET: Parametrization of Cosserat Equations, Acta Mechanica, 215, 2010, 43-55.\\
\noindent
[21] J.-F. POMMARET: Macaulay Inverse Systems Revisited, Journal of Symbolic Computation, 46, 2011, 1049-1069.\\
\noindent
[22] J.-F. POMMARET: Spencer Operator and Applications: From Continuum Mechanics to Mathematical Physics, in "Continuum Mechanics-Progress in Fundamentals and Engineering Applications", Dr. Yong Gan (Ed.), ISBN: 978-953-51-0447--6, InTech, 2012, Available from: \\
http://www.intechopen.com/books/continuum-mechanics-progress-in-fundamentals-and-engineerin-applications/spencer-operator-and-applications-from-continuum-mechanics-to-mathematical-physics  \\
\noindent
[23] J.-F. POMMARET: A Pedestrian Approach to Cosserat/Maxwell/Weyl Theory, Preprint 2012 \\
http://hal.archives-ouvertes.fr/hal-00740314     \\
http://fr.arXiv.org/abs/1210.2675           \\
\noindent
[24] J.-F. POMMARET, A. QUADRAT: Localization and parametrization of linear multidimensional control systems, Systems \& Control Letters 37 (1999) 247-260.  \\
\noindent
[25] J.-F. POMMARET, A. QUADRAT: Algebraic Analysis of Linear Multidimensional  Control Systems, IMA Journal of Mathematical Control and Informations, 16, 1999, 275-297.\\
\noindent
[26] A. QUADRAT: Analyse Alg\'ebrique des Syst\`emes de Contr\^ole Lin\'eaires Multidimensionnels, Th\`ese de Docteur de l'Ecole Nationale 
des Ponts et Chauss\'ees, 1999 \\
(http://www-sop.inria.fr/cafe/Alban.Quadrat/index.html).\\
\noindent
[27] A. QUADRAT: Une Introduction \`{a} l'Analyse Alg\'{e}brique Constructive et \`{a} ses Applications, INRIA Research Report 7354, AT-SOP Project, july 2010. Les Cours du CIRM, 1 no. 2: Journ\'{e}es Nationales de Calcul Formel (2010), p281-471 (doi:10.5802/ccirm.11). \\
\noindent
[28] J.J. ROTMAN: An Introduction to Homological Algebra, Pure and Applied Mathematics, Academic Press, 1979.\\
\noindent
[29] W. M. SEILER: Involution: The Formal Theory of Differential Equations and its Applications to Computer Algebra, Springer, 2009, 660 pp. (See also doi:10.3842/SIGMA.2009.092 for a recent presentation of involution, in particular sections 3 (p 6 and reference [11], [22]) and 4). \\
\noindent
[30] D.C. SPENCER: Overdetermined Systems of Partial Differential Equations, Bull. Amer. Math. Soc., 75, 1965, 1-114.\\
\noindent
[31] O. ZARISKI, P. SAMUEL: Commutative Algebra, Van Nostrand, 1958.  \\
\noindent
[32] E. ZERZ: Topics in Multidimensional Linear Systems Theory, Lecture Notes in Control and Information Sciences 256, Springer, 2000.\\

\end{document}